\documentclass[12pt,a4paper]{article}
\usepackage{amsmath,amsfonts,amssymb,amsthm,amscd}
\usepackage{hyperref}
\usepackage{braket}
\usepackage{epstopdf}

\usepackage{braket}
\usepackage{xypic}
\usepackage{amsmath}

\usepackage[american]{babel}
\usepackage{color} 
\usepackage{epsfig}
\usepackage{caption}
\usepackage{rotating}
\usepackage{setspace}
\usepackage{fancyhdr}
\usepackage{booktabs}
\usepackage{mathrsfs}
\usepackage{amsthm}
\usepackage{amsmath,amssymb}
\usepackage{enumerate}
\usepackage{tikz}
\usepackage{youngtab}
\usepackage{hyperref}
\usepackage{tensor}
\usepackage{mathabx}
\usepackage{booktabs}
\usepackage{mathrsfs}
\usepackage{subfigure}
\usepackage{amsmath,amssymb}

\def\id{{\rm id}}

\renewcommand{\[}{\begin{equation}}
\renewcommand{\]}{\end{equation}}

\newtheorem{thm}{Theorem}[section]
\newtheorem{cor}[thm]{Corollary}
\newtheorem{lem}[thm]{Lemma}
\newtheorem{prop}[thm]{Proposition}

\theoremstyle{definition}
\newtheorem{defn}{Definition}[section]

\theoremstyle{remark}
\newtheorem{oss}{Remark}[section]

\newcommand{\ra}{\rightarrow}

\DeclareMathOperator{\supp}{supp}
\DeclareMathOperator{\SL}{SL}

\DeclareMathOperator{\spann}{span}

\DeclareMathOperator{\Q}{Q-}

\DeclareMathOperator{\cohull}{cohull}

\title{\huge An approach to the study of boundary actions}

\author{J. Bassi\\ \\Department of Mathematics, University of Tor Vergata, \\Via della Ricerca Scientifica 1, 00133, Roma, Italy\\ \\email: bssjcp01@uniroma2.it}

\begin{document}

\date{}
\maketitle

\begin{abstract}
Given an action of a discrete countable group $G$ on a countable set $\mathfrak{X}$, it is studied the relationship between properties of the associated Calkin representation and the dynamics of the group action on the boundary of the Stone-\v{C}ech compactification of $\mathfrak{X}$. The first section contains results about amenability properties of actions of discrete countable groups on non-separable spaces and is of independent interest. In the second section these results are applied in order to translate regularity properties of the Calkin representation and the topological amenability on the Stone-{\v C}ech boundary within the common framework of measurable dynamics on certain extensions of the Stone-{\v C}ech boundary of $\mathfrak{X}$. 
\end{abstract}

\tableofcontents




\section*{Introduction}

Regularity properties of boundary actions are of fundamental importance in the study of discrete groups; in particular amenability in the sense of R. J. Zimmer (\cite{Zi}) and (topological) amenability in the sense of C. Anantharaman-Delaroche (\cite{An}) turn out to play a crucial role in this setting. Among the first evidences of the importance of these concepts it is worth mentioning the celebrated work \cite{adams}, where S. Adams proved that any quasi-invariant finite measure on the Gromov boundary of an hyperbolic group is Zimmer-amenable and deduced rigidity results from this, and the manuscript \cite{Oz1} by N. Ozawa, in which it is shown that exact groups are precisely the groups for which the left action on the boundary of their Stone-{\v C}ech compactification is topologically amenable. A crucial role in the recent development of boundary theory is played by the work \cite{KaKe}, where  M. Kalantar and M. Kennedy proved that the Furstenberg boundary of a discrete countable group is the spectrum of the equivariant injective envelope of the complex numbers and used this fact in order to establish many interesting fact, for example they answered in the positive a conjecture made by N. Ozawa in \cite{Oz2} (see also \cite{BaRa1} for a related topic). This approach led to impressive results in \cite{BrKaKeOz}.\\Topological amenability on Stone-{\v C}ech boundaries for (non necessarily free and transitive) actions of discrete countable groups on countable sets is also of great impact; for example it was used in \cite{Oz} in order to obtain prime factorizations results, in \cite{BaRa2} in order to confirm a prediction by B. Bekka and M. Kalantar from \cite{BeKa}, and in \cite{BaRa} in order to study a mixing property associated to a certain action of $\SL(3,\mathbb{Z}) \times \SL(3,\mathbb{Z})$. A particular case of interest regards the left-right action of a discrete countable group on its Stone-{\v C}ech boundary; the topological amenability of this action is usually referred to as {\it biexactness} or {\it property $\mathcal{S}$}. This property has been proved in the case of hyperbolic groups (and certain variations on this theme) (see for example \cite{BrOz}) and it implies the {\it Akemann-Ostrand property} (\cite{AkOs}), i.e. the temperedness of the canonical representation of the full group $C^*$-algebra of the cartesian product of the group with itself into the Calkin algebra of the group. The latter property reflects a very strong property at the level of the group von Neumann-algebra, i.e. {\it solidity}. In particular, the von Neuman algebras associated to groups satisfying the Akemann-Osrtand property do not lie in the realm of crossed product von Neuman algebras (note anyway that in some specific cases they can be reconstructed from certain associated crossed products, as observed in \cite{Ba1}, using the results in \cite{Ba} and \cite{GoLi}). If a property $T$ group satisfies the Akemann-Ostrand property, then it is not nuclear in $K$-theory (\cite{Sk} Corollary 4.5). Weakly amenable, non-amenable biexact groups are $\mathcal{C}$-rigid (\cite{PoVa}). A weakening of biexactness was considered in \cite{BIP}, where some rigidity results have been obtained for a much larger class of groups. The question whether biexactness and the Akemann-Ostrand property are equivalent is still open and there are neither evidences nor partial results towards a possible answer. \\This work originates from the will to find a solution to this problem and aims to a better understanding of the relationship between regularity properties of the Calkin representation associated to the action of a discrete countable group on a countable set and dynamical properties of the action of the group on the Stone-{\v C}ech boundary of this set. This manuscript represents the first step in this direction and contains the general framework which will be used to study specific cases in the future.

Section \ref{am} contains results about amenable actions of discrete countable groups on compact Hausdorff spaces and is of independent interest. The characterization given in this section of topological amenability is the starting point for the considerations appearing in the subsequent sections. Namely, it is proved that the action of a discrete countable group $G$ on a compact Hausdorff space $X$ is topologically amenable if and only if for every quasi-invariant probability measure $\mu$ on $X$ the following two conditions are satisfied:
\begin{itemize}
\item[(i)] Pair amenability: there is an equivariant projection $L^\infty (\mu) \overline{\otimes} L^\infty (\mu) \rightarrow L^\infty (\mu)$,
\item[(ii)] Temperedness: the Koopman representation associated to $\mu$ is tempered.
\end{itemize}
This result is in complete analogy with the results contained in \cite{BEW} and \cite{An}.

In Section \ref{sec2} we use the above result in order to find a common framework for the study of the Calkin representation associated to the action of $G$ on a countable set $\mathfrak{X}$ and topological amenability of the action of $G$ on $\partial_\beta \mathfrak{X}$. In order to do so, we use in Section \ref{susec2.3} considerations from \cite{BaRa} in order to exhibit a particular set of quasi-invariant measures defined on certain "non-standard" extensions of the Stone-{\v C}ech boundary of $\mathfrak{X}$ indexed by free ultrafilters over $\mathbb{N}$ and denoted by $\partial_{\beta, \omega} \mathfrak{X}$, with $\omega \in \partial_\beta \mathbb{N}$. The temperedness of all these measures corresponds to the temperedness of the associated Calkin representation $C^* (G) \rightarrow \mathbb{B}(l^2(\mathfrak{X}))/\mathbb{K}(l^2(\mathfrak{X}))$. Moreover, at the beginning of Section \ref{sec2} it is shown that topological amenability of the action of $G$ on $\partial_\beta \mathfrak{X}$ is equivalent to topological amenability of the action of $G$ on any of these extensions. Hence  we conclude that the two concepts can be described in the common framework of quasi-invariant measures on the aforementioned extensions of $\partial_\beta \mathfrak{X}$. Furthermore, for all of these extensions we exhibit a small  set of quasi-invariant probability measures for which pair amenabiity is automatic (\ref{sec2.1}).

\section*{Notation}
\label{notation} If $A$ is a $C^*$-algebra we denote by $\sigma (A)$ its spectrum. If $\mu$ is a Radon measure on a locally compact Hausdorff space $Z$, we denote by $\eta_\mu : C_0(Z) \rightarrow L^2(\mu)$ the canonical linear map sending a continuous function to its $L^2(\mu)$-equivalence class and by $\pi_\mu : C(X) \rightarrow \mathbb{B}(L^2 (\mu))$ the $*$-representation given by multiplication; we also denote by $\braket{\cdot , \cdot}_\mu$ the associated scalar product in $L^2(\mu)$ and by $\| \cdot \|_{2, \mu}$ the associated Hilbert space norm. We denote by $\mathcal{P}(X)$ the set of Radon probability measures on $X$ and if $X$ is endowed with the action of a group $G$, we denote by $\Q\mathcal{P}(X) \subset \mathcal{P}(X)$ the set of quasi-invariant probability measures.

\section{Amenability for non-separable dynamical systems}
\label{am}
\subsection{Amenability and pair amenability}
Let $X$ be a locally compact Hausdorff space and $G$ a discrete countable group acting on $X$. 
\begin{defn}[\cite{AD2}, Definition 2.1]
The action of $G$ on $X$ is said \textit{topologically amenable} if there is a net of continuous functions $\mu_\lambda : X \rightarrow \mathcal{P}(G)$ such that 
\begin{equation*}
\lim_\lambda \sup_{x \in Y} \| s \mu (x) - \mu (s x)\|_1 =0
\end{equation*}
for every compact set $Y \subset X$ and $s \in G$.
\end{defn}
Suppose now that $\mu$ is a positive Radon measure on $X$ which is quasi-invariant with respect to the action of $G$. Then there is a well defined action on $L^\infty (X,\mu)$.
\begin{defn}[\cite{AD3}, D{\'e}finition 3.4]
The action of $G$ on $L^\infty (X,\mu)$ is \textit{Zimmer-amenable} if there is an equivariant norm-$1$ projection
\begin{equation*}
P: l^\infty (G) \overline{\otimes}L^\infty (X,\mu) \rightarrow L^\infty (X,\mu),
\end{equation*}
where the von Neumann tensor product $l^\infty (G) \overline{\otimes}L^\infty (X,\mu)$ is endowed with the diagonal action of $G$. If this is the case we will say that the measure $\mu$ is Zimmer-amenable.
\end{defn}
In the above situation the von Neumann algebra $L^\infty (X,\mu)$ is represented on $L^2(X,\mu)$ by multiplication and is in fact the weak closure of the image of $C_0 (X)$ under this representation. The representation by multiplication for $C_0(X)$ actually gives rise to a covariant representation.
\begin{defn}
\label{koop}
Let $\mu$ be a $G$-quasi-invariant positive Radon measure on $X$. The \textit{Koopman representation} associated to $\mu$ is the group homomorphism $u^\mu:G \rightarrow U(L^2(\mu))$ given by
\begin{equation*}
u^\mu_s \xi = \sqrt{\frac{ds\mu}{d\mu}} \cdot \xi \circ s^{-1} \qquad \mbox{ for } \xi \in L^2(\mu).
\end{equation*}
\end{defn}

It is well known that if the action of $G$ on $X$ is topologically amenable, then every $G$-quasi-invariant positive Radon measure on $X$ is Zimmer-amenable (see for example \cite{BEW} Lemma 3.28). Moreover, if $\mu$ is Zimmer-amenable, then $u^\mu$ is tempered (see for example \cite{BaRa} Lemma 2.2). The problems of determining under which conditions any of these implications can be reversed is highly non-trivial. Hence it is natural to introduce the following 
\begin{defn}[cf. \cite{BEW} Definition 3.25]
The action of $G$ on $X$ is said \textit{measurewise amenable} if for every $G$-quasi-invariant Radon probability measure $\mu$ on $X$ it holds that $\mu$ is Zimmer-amenable.
\end{defn}

The problem of the equivalence between topological amenability and measurewise amenability has been addressed in \cite{BEW}, where a positive answer has been given in the case of a second countable group acting on a second countable space (\cite{BEW} Corollary 3.29). In Corollary \ref{cor1.6} it is proved an analogue result for the case of a discrete countable group acting on a (in general non-separable) compact Hausdorff space, which will be used in the forthcoming sections. 

On the other hand the question of whether the temperedness of the Koopman representations associated to $G$-quasi-invariant measures implies the Zimmer-amenability of these measures has a negative answer in general. This problem was studied in \cite{AD1} again in the second countable case. The approach there was to translate these forms of regularity in properties of groupoid representations. An important ingredient is the following
\begin{defn}
Let $M \subset N$ be an inclusion of von Neumann algebras and suppose $N$ is endowed with an action of a discrete countable group $G$ which restrict to an action on $M$. We say that the pair $(N,M)$ is \textit{amenable} if there is an equivariant projection $E : N \rightarrow M$.
\end{defn}
In \cite{AD1} Theorem 4.3.1 it is proved that for actions of second countable locally compact groups acting on standard Borel spaces, Zimmer-amenability of the quasi-invariant measure $\mu$ is equivalent to the fact that the Koopman representation is tempered and the pair $(L^\infty (X,\mu) \overline{\otimes} L^\infty (X,\mu), L^\infty (X,\mu))$ is amenable. The analogue in our situation is Corollary \ref{cor1.4}.

\subsection{A characterization of Zimmer-amenability}

Let $X$ be a locally compact Hausdorff space, $\mu$ a positive Radon measure on $X$ and $E$ be a Banach space. A map $f : X \rightarrow E$ is said to be $\mu$-almost surely separately valued if there is a measurable set $ A \subset X$ such that $\mu (X \backslash A)=0$ and the range of $f|_A$ is separable in $E$. We say that $f$ is $\mu$-measurable if it is $\mu$-almost surely separately valued and $f^{-1} (U)$ is measurable for every open set $U \subset E$. We refer to \cite{Ta1} IV.7 and \cite{Wil} B.1 for the theory of $\mu$-measurable functions. The following properties hold:
\begin{itemize}
\item[(P1)] If $f : X \rightarrow M$ is $\mu$-measurable, the same is true for $\|f\| : x \mapsto \|f(x)\|$ (\cite{Wil} Lemma B.29),
\item[(P2)] if $E$ is a Banach algebra and $f,g$ are $\mu$-measurable, then the pointwise product $fg$ is $\mu$-measurable (\cite{Wil} Corollary B.21)
\end{itemize}
In virtue of (P1) it makes sense to define the spaces 
\begin{equation*}
L^1((X,\mu),E) := \{ f: X \rightarrow E \; | \; f \mbox{ is $\mu$-measurable and } \int_X \|f(x)\| d\mu  < \infty \}/\sim ,
\end{equation*}
\begin{equation*}
L^2((X,\mu),E) := \{ f: X \rightarrow E \; | \; f \mbox{ is $\mu$-measurable and } \int_X \|f(x) \|^2 d\mu  < \infty \}/\sim ,
\end{equation*}
where in both cases "$\sim$" is equivalence $\mu$-almost everywhere. The spaces $L^1((X,\mu),E)$ and $L^2((X,\mu),E)$ are complete (\cite{Ta1} Proposition IV.7.4) and they share the following properties:
\begin{itemize}
\item[(P3)] the set $C_c (X,E)$ of continuous compactly supported $E$-valued functions is dense in both $L^1((X,\mu),E)$ and $L^2((X,\mu),E)$ (indeed these Banach spaces are defined as appropriate completions of $C_c (X,E)$ in \cite{Ta1} Section IV.7),
\item[(P4)] the simple $E$-valued functions are dense in both $L^1((X,\mu),E)$ and $L^2((X,\mu),E)$ (the statement for $L^1((X,\mu), E)$ is proved in \cite{Wil} Proposition B.33, but a similar argument can be applied to $L^2((X,\mu), E)$),
\item[(P5)] there is a contractive linear map $L^1 ((X,\mu),E) \rightarrow E$, $f \mapsto \int_X f(x) d\mu (x)$, which is characterized by the following property: for every $f \in L^1((X,\mu),E)$
\begin{equation*}
\phi (\int_X f(x) d\mu (x)) = \int_X \phi (f(x)) d\mu (x) \qquad \mbox{ for all } \phi \in E^*
\end{equation*}
(\cite{Wil} Proposition B.34).
\end{itemize}

We will be interested in the case in which $E=M$ is a von Neumann algebra and stick to this case in the following. It follows from property (P2) and the fact that products of $\mu$-measurable functions respect equivalence $\mu$-almost everywhere that it make sense to multiply pointwise elements in $L^2((X,\mu),M)$, the result being an element in $L^1((X,\mu),M)$. For $f,g \in L^2((X,\mu),M)$, we denote by $\langle f,g \rangle_M$ the $M$-valued scalar product given by $\langle f,g \rangle = \int_X f^* g d\mu$. \\Let $L^1((X,\mu)M)_1$ be the set of functions in $L^1((X,\mu),M)$ such that $\int_X f(f) d\mu (x) =1 \in M$ and let $L^1((X,\mu),M)_{1, +}$ be the subset of those functions for which such representatives can be chosen pointwise positive. Similarly, denote by$L^2((X,\mu),M)_{1}$ the set of functions $f$ in $L^2((X,\mu),M)$ such that $\langle f, f \rangle_M=1 \in M$ and by $L^2((X,\mu),M)_{1,+}$ the subset of functions in $L^2((X,\mu),M)$ admitting representatives which are pointwise positive. Given an element $f \in L^{p}((X,\mu),M)$, $p\in \{1,2\}$ we will denote by $\|f\|_{p,\infty}$ its norm $(\int_X \|f(x)\|^p d\mu (x))^{1/p}$; if the meaning is clear from the context, we will just write $\|f\|$.\\


Let $M$ be a commutative von Neumann algebra. We will denote by $M_*$ its predual and by $M_{*,+}$ its positive part. In virtue of \cite{Ta1} Theorem IV.7.17 and Lemma 7.13 for every $F \in L^\infty (X,\mu) \bar{\otimes} M = L^1((X,\mu), M_*)^*$ there is a uniformly bounded function $\tilde{F} : X \rightarrow M$ such that for every $\omega \in M_*$ the function $X \rightarrow \mathbb{C}$, $x \mapsto \langle \tilde{F} (x), \omega \rangle$ is measurable and $\langle F, g \rangle = \int_X \langle \tilde{F} (x), g (x) \rangle$ for every $g \in L^1((X,\mu),M_*)$, moreover $\|\tilde{F}\|_\infty \leq \| F\|$ (note that we cannot assume $\tilde{F}$ to be chosen almost surely separately valued, but it will satisfy the property that the map $x \mapsto \langle \tilde{F} (x), \omega \rangle$ be measurable for every $\omega \in M_*$, together with the measurability property of \cite{Ta1} Definition IV.7.7). 
Let now $h \in L^1 ((X,\mu),M)_{1,+}$; we proceed as in \cite{BeCr} 3.1, hence we want to associate to $h$ a unital positive $M$-bimodule map $P_h : L^\infty (X,\mu) \bar{\otimes} M \rightarrow M$. Given $\omega \in M_*$ let $L^1((X,\mu),M) \rightarrow L^1((X,\mu), M_*)$, $g \mapsto \tilde{\omega}(g)$ be the linear map given in \cite{BeCr} Lemma 2.1. We define $P_h : L^\infty (X,\mu) \bar{\otimes} M \rightarrow M$, $F \mapsto \{ \omega \mapsto \langle F, \tilde{\omega} (h)\rangle\}$; explicitely, this is given by 
\begin{equation*}
\begin{split}
\langle P_h (F), \omega \rangle &= \int_X \langle \tilde{F}(x) , \tilde{\omega}(h) (x) \rangle d\mu (x)\\ &= \int_X \langle \tilde{F}(x) h(x), \omega \rangle d\mu (x), \qquad F \in L^\infty (X) \bar{\otimes} M, \; \omega \in M_* .
\end{split}
\end{equation*}
Since the above result does not depend on the choice of the representative $\tilde{F}$ for $F$, we will unambiguously write $\langle P_h (F), \omega \rangle = \int_X \langle F(x) h(x), \omega \rangle d\mu (x)$. Since $P_h$ is a positive contractive $M$-bimodule map, it is completely positive.

Denote by $\mathcal{B}(L^\infty(X,\mu) \overline{\otimes} M, M)$ the set of bounded $M$-bimodule maps from $L^\infty (X,\mu) \overline{\otimes} M$ to $M$, by $\mathcal{P}$ the subset of positive $M$-bimodule maps of norm $\leq 1$ and by $\mathcal{P}_1$ the subset of $\mathcal{P}$ of maps given by elements in $L^1((X,\mu),M)_{1,+}$ following the above procedure.

\begin{lem}
\label{lem3.3}
The convex set $\mathcal{P}_1$ is dense in $\mathcal{P}$ for the point-weak$^*$ topology.
\end{lem}
\proof The proof of the analogous result for the set $\mathcal{P}_K$ in \cite{BeCr} Lemma 3.1 uses the fact that the bipolar of $\mathcal{P}_K$ contains $\mathcal{P}$; the result follows since $\mathcal{P}_K \subset \mathcal{P}_1$. $\Box$

In the following we will be interested in the case in which a von Neumann algebra $M$ is endowed with the action of a discrete countable group $G$ and the same group also acts on the space $X$ leaving invariant the measure class of $\mu$, i.e. $\mu$ is $G$-quasi-invariant. Given $s \in G$, we will write $a \mapsto sa$ for $a$ either an element of $M$ or of $X$ (hence also of $\mathcal{P}(X)$) for the corresponding translates under the group action. We will also write $s \otimes s$ for the tensor product actions of $G$. Note that if $\mu$ is $G$-quasi-invariant then $s\mu \precsim \mu$ for every $s \in G$ and so it makes sense to consider the Radon-Nykodim derivative $\sqrt{ds\mu /d\mu} \in L^1(X,\mu)$.

\begin{oss}
Given a von Neumann algebra $M$ (commutative or not) endowed with the action of a discrete countable group $G$, it is possible to define a dual action on its dual $M^*$ given by $(s,\omega) \mapsto \omega \circ s$. Let now $X$ be a locally compact Hausdorff space, $\mu$ a (possibly infinite) positive $G$-quasi-invariant Radon measure on $X$, $M$ a commutative von Neumann algebra endowed with an action of $G$ and let $g \in L^1((X,\mu), M_*)$. For every $ f \in L^\infty (X,\mu)$, $m \in M$ and $s \in G$ we have 
\begin{equation*}
\begin{split}
\int_X &(s\otimes s) (f \otimes m) g d\mu (x)= \int_X f(s^{-1} x) \langle sm, g(x)\rangle d\mu (x) \\ &= \int_X \frac{ds^{-1} \mu}{d\mu} (x) f(x) \langle sm, g(sx)\rangle d\mu (x) = \int_X  f(x) \langle sm, \frac{ds^{-1} \mu}{d\mu} (x) g(sx)\rangle d\mu (x); 
\end{split}
\end{equation*}
note that the function $x \mapsto \frac{ds^{-1} \mu}{d\mu} (x) (g(sx) \circ s)$ is $\mu$-measurable in virtue of property (P2) above and has the same $L^1$-norm of the original function $g$. Since the linear span of elements of the form $f \otimes m$ is dense in $L^\infty (X,\mu) \overline{\otimes} M$, it follows that $g \circ (s \otimes s) =  \{ x \mapsto \frac{ds^{-1} \mu}{d\mu} (x) (g(sx) \circ s)\}$. It follows from this that if $F \in L^\infty (X,\mu) \overline{\otimes} M$ and $\tilde{F}$ is a measurable representative, then for every $g \in L^1((X,\mu), M_*)$ we have $\langle(s \otimes s) F, g \rangle = \langle F, g \circ (s \otimes s) \rangle = \int_X \langle \tilde{F} (x),  \frac{ds^{-1} \mu}{d\mu} (x) (g(sx) \circ s)\rangle d\mu (x) = \int_X \langle s(\tilde{F}(s^{-1} x)), g(x)\rangle d\mu (x)$, thus $\{ x \mapsto s(\tilde{F}(s^{-1} x))\}$ is a representative for $(s \otimes s) F$ (which again satisfies the property that $x \mapsto \langle s(\tilde{F}(s^{-1} x)), \omega \rangle $ is measurable for every $\omega \in M_*$ and the measurability assumption of \cite{Ta1} Definition IV.7.7).
\end{oss}

\begin{prop}
\label{prop3.4}
Let $M$ be a commutative von Neumann algebra, $X$ be locally compact Hausdorff space and $G$ a discrete countable group. Let $G$ act on both $M$ and $X$ via $(s,m) \mapsto sm$ and $(s,x) \mapsto sx$ respectively . Let $\mu$ be a $G$-quasi-invariant (finite or infinite) positive Radon measure on $X$. Endow $L^\infty (X,\mu) \overline{\otimes}M$ with the diagonal action of $G$. The following are equivalent:
\begin{itemize}
\item[(i)] there is a $G$-equivariant norm one projection from $L^\infty (X,\mu) \overline{\otimes} M$ to $1 \overline{\otimes} M$,
\item[(ii)] there is a net $(h_\lambda)$ in $L^1 ((X,\mu),M)_{1, +}$ such that 
\begin{equation*}
\int_X | \frac{d s\mu}{d\mu} s (h_\lambda (s^{-1} x)) - h_\lambda (x)| d\mu(x) \rightarrow 0 \qquad \mbox{weak}^* \quad \forall s \in G,
\end{equation*}
\item[(iii)] there is a net $(\xi_\lambda) \in L^2((X,\mu), M)_{1,+}$ such that
\begin{equation*}
\langle \xi_\lambda, \sqrt{\frac{ds\mu }{d\mu}} s(  \xi_\lambda (s^{-1} \cdot))\rangle_M  \rightarrow 1 \qquad \mbox{weak}^* \quad \forall s \in G,
\end{equation*}
\item[(iv)]  there is a net $(\xi_\lambda) \in L^2((X,\mu), M)_{1}$ such that
\begin{equation*}
\langle \xi_\lambda , \xi_\lambda\rangle_M  \rightarrow 1 \qquad \mbox{ weak}^* \qquad \mbox{ and}
\end{equation*}
\begin{equation*}
\langle \xi_\lambda, \sqrt{\frac{ds\mu }{d\mu}} s( \xi_\lambda (s^{-1} \cdot))\rangle_M  \rightarrow 1 \qquad \mbox{weak}^* \quad \forall s \in G.
\end{equation*}

\end{itemize}

\end{prop}
\label{prop2.10}
\proof (i)$\Rightarrow$(ii): Let $P$ be an equivariant norm one projection from $L^\infty (X)\bar{\otimes}M$ to $M$ and consider a net $(g_\lambda)$ in $L^1((X,\mu),M)_{1,+}$ realizing the approximants as in Lemma \ref{lem3.3}. For every $F \in L^\infty (X) \overline{\otimes} M$ and $\omega \in M_*$ we have
\begin{equation*}
\braket{P (F), \omega} = \lim_\lambda \int_X \braket{F(x) g_\lambda (x), \omega} d\mu (x).
\end{equation*}
The equivariance of $P$ means that for every $s \in G$ we have
\begin{equation*}
\lim_\lambda \int_X \braket{g_\lambda (x) s(F(s^{-1} x)), \omega} d\mu (x) = \lim_\lambda  \int_X \braket{g_\lambda (x) F( x),   \omega \circ s }
\end{equation*}
On the other hand we have
\begin{equation*}
\begin{split}
\int_X \braket{g_\lambda (x) s(F(s^{-1} x)), \omega} d\mu (x)= \int_X \frac{d s^{-1} \mu}{d\mu} (x) \braket{(g_\lambda (s x) s(F(x)), \omega} d \mu (x)\\
= \int_X \frac{d s^{-1} \mu}{d\mu} (x) \braket{(s^{-1}g_\lambda (s x)) F(x),  \omega \circ s} 
\end{split}
\end{equation*}
Thus 
\begin{equation*}
\frac{ds\mu}{d\mu} s (g_\lambda (s^{-1} \cdot )) - g_\lambda (\cdot ) \rightarrow 0 \qquad \mbox{ for the topology } \tau_f, \; \forall s \in G
\end{equation*}
where $\tau_f$ is the topology on $L^1((X,\mu),M)$ induced by the family of seminorms of the form $p_{F,\omega} (g) = | \int_X \braket{g(x) F(x), \omega} d\mu (x)|$, with $F \in L^\infty (X,\mu) \overline{\otimes} M$, $\omega \in M_*$.  Let $F \subset G$ and $A \subset M_{*,+}$ be finite sets. The net $\oplus_{s \in F} (\frac{ds\mu}{d\mu} s (g_\lambda (s^{-1} \cdot )) - g_\lambda (\cdot ))$ in $\oplus_{s \in F} L^1 ((X,\mu), M)$ converges to zero with respect to the topology $\tau_f^{|F|}$. For $\omega \in M_{*,+}$ let $\tilde{\omega} : L^1((X,\mu),M) \rightarrow L^1((X,\mu),M_*)$ be the continuous linear map given in \cite{BeCr} Lemma 2.1. By linearity the image under $\oplus_{\omega \in A} \tilde{\omega}^{|F|}$ of the convex hull of $\{\oplus_{s \in F} (\frac{ds\mu}{d\mu} s (g_\lambda (s^{-1} \cdot )) - g_\lambda (\cdot ))\}$ is convex. The fact that the net converges to zero in $\tau_f^{|F|}$ implies that the image of the net under $\oplus_{\omega \in A} \tilde{\omega}^{|F|}$ converges to zero weakly (since $L^1(X, M_*)^* = L^\infty (X) \overline{\otimes} M$ and the dual of a finite direct sum is the direct sum of the duals) and so the image of the convex hull of the set $\{\oplus_{s \in F} (\frac{ds\mu}{d\mu} s (g_\lambda (s^{-1} \cdot )) - g_\lambda (\cdot ))\}$ under $\tilde{\omega}$ contains zero in its weak closure, hence contains zero in its norm closure. Thus given $\epsilon >0$ we can find $h$ in the convex hull of $\{g_\lambda\}$ such that $\sum_{\omega \in A} \sum_{s \in F} \int \| \tilde{\omega} (\frac{ds\mu}{d\mu} s (h \circ s^{-1} )- h)(x)\| <\epsilon$. Now, for every $h \in L^1((X,\mu),M)$ we have for $\mu$-almost every $x \in X$ that $\|\tilde{\omega}(h)(x)\| = \sup_{m \in M_1} |\omega (h(x)m)| \geq \omega (h(x) u (x))= \omega (|h(x)|)$, where $u(x)$ is a partial isometry in $M$ satisfying $h(x) u(x) = |h(x)|$. We thus have
\begin{equation*}
\int_X \omega (|\frac{ds\mu}{d\mu} s (h \circ s^{-1}(x) )- h (x)|) d\mu (x)\leq \int_X \|\tilde{\omega} (\frac{ds\mu}{d\mu} s (h \circ s^{-1} )- h)(x)\| d\mu (x)< \epsilon 
\end{equation*}
for every $\omega \in A$, $s \in F$. The result follows from (P5).\\
(ii)$\Rightarrow$(iii): Let $(h_\lambda)$ be a net as in the statement and define $\xi_\lambda (x) = \sqrt{h_\lambda (x)}$. This net takes values in $L^2((X,\mu),M)_{1, +}$. For $\omega \in M_{*,+}$ we have
\begin{equation*}
\begin{split}
|\omega (\braket{\xi_\lambda, \sqrt{\frac{ds\mu}{d\mu}} s(\xi_\lambda \circ s^{-1} )} - 1))|&= |\int_X \omega (\xi_\lambda(x) [\sqrt{\frac{ds\mu}{d\mu}} (x)s(\xi_\lambda (s^{-1}x)) -\xi_\lambda (x)] d\mu (x)|\\&\leq (2\|\omega\| \int_X \omega(|\frac{ds\mu}{d\mu} (x) s(h_\lambda (s^{-1}x)) - h_\lambda (x)|)d\mu(x))^{1/2}
\end{split}
\end{equation*}
by Cauchy-Shwartz inequality; the result follows.\\
(iii)$\Rightarrow$(ii): Let $(\xi_\lambda)$ be a net in $L^2((X,\mu),M)_{1,+}$ as in the statement and $\omega \in M_{*,+}$, $s \in G$. We have
\begin{equation*}
\begin{split}
\int_X \omega&(| \frac{ds\mu}{d\mu}(x) s(\xi_\lambda^2 (s^{-1}x)) - \xi_\lambda^2 (x)|)d\mu (x) \\ &= \int_X \omega(| \sqrt{\frac{ds\mu}{d\mu} (x)}s(\xi_\lambda (s^{-1}x)) + \xi_\lambda (x)| | \sqrt{\frac{ds\mu}{d\mu}(x) }s(\xi_\lambda (s^{-1}x))-\xi_\lambda (x)| ) d\mu(x)\\
&\leq (4 \|\omega\| \int_X \omega(| \sqrt{\frac{ds\mu}{d\mu}(x) }s(\xi_\lambda (s^{-1}x))-\xi_\lambda (x)|^2 ) d\mu(x))^{1/2}\\
&= (4 \|\omega\|[\int_X \omega (\frac{ds\mu}{d\mu}(x) s(\xi_\lambda^2 (s^{-1}x)) d\mu(x)+\int_X \omega(\xi_\lambda^2 (x))d\mu(x) \\&- 2\int_X \omega ( \xi_\lambda (x)\sqrt{\frac{ds\mu}{d\mu}(x) }s(\xi_\lambda (s^{-1}x)))d\mu(x)])^{1/2} \rightarrow 0.
\end{split}
\end{equation*}
(ii)$\Rightarrow$(i): By compactness of $\mathcal{P}$ in the point-weak$^*$ topology the net $(P_{g_\lambda})$ admits a convergent subnet, which is a projection of norm $1$. After reindexing we can suppose $P_{g_{\lambda}} \rightarrow P$. We need to check that $P$ is equivariant. For let $s \in G$, $\omega \in M_{*,+}$, $F \in L^\infty (X,\mu) \overline{\otimes} M$ and compute
\begin{equation*}
\begin{split}
|\lim_\lambda &[ \int_X \omega (g_\lambda (x)s(F (s^{-1}x))) d\mu(x) - \int_X (\omega \circ s) (g_\lambda (x) F(x)) d\mu(x)]|\\
&= |\lim_\lambda [\int_X (\omega \circ s) (s^{-1} (g_\lambda (x)) F(s^{-1} x) d\mu(x) -\int_X (\omega \circ s) (g_\lambda (x) F(x)) d\mu(x)]|\\
& = | \lim_\lambda \int_X (\omega \circ s) ((\frac{ds^{-1}\mu}{d\mu}(x) s^{-1} (g_\lambda (sx)) -g_\lambda (x)) F(x)) d\mu(x)|\\
&\leq \|\tilde{F}\|_\infty \lim_\lambda \int_X |(\omega \circ s) (\frac{ds^{-1}\mu}{d\mu}(x) s^{-1} (g_\lambda (sx)) -g_\lambda (x))| d\mu(x)\\
& \leq \|\tilde{F}\|_\infty \lim_\lambda \int_X (\omega \circ s) (|\frac{ds^{-1}\mu}{d\mu}(x) s^{-1} (g_\lambda (sx)) -g_\lambda (x)|) d\mu(x)=0. 
\end{split}
\end{equation*}
(iv)$\Rightarrow$(ii): Let $(\xi_\lambda)$ as in the hypothesis and consider $\eta_\lambda =\xi_\lambda^* \xi_\lambda$ for every $\lambda$. Let $\omega \in M_{*,+}$. We have $\omega(\int_X \eta_\lambda (x) d\mu (x)) \rightarrow 1$ and
\begin{equation*}
\begin{split}
\int_X \omega &(|\frac{ds\mu}{d\mu} (x) s(\xi_\lambda^* (s^{-1} x)) s(\xi_\lambda (s^{-1} x)) - \xi_\lambda^* (x) \xi_\lambda (x)|)d\mu (x)\\
&=  \int_X \omega (|\sqrt{\frac{ds\mu}{d\mu} (x)} s(\xi_\lambda^* (s^{-1} x)) - \xi_\lambda^* (x) | |\sqrt{\frac{ds\mu}{d\mu} (x)} s(\xi_\lambda (s^{-1}x)) + \xi_\lambda (x)|) d\mu(x)\\
&\leq (4 \|\omega \| \|\xi_\lambda\|^2 \int_X \omega (|\sqrt{\frac{ds\mu}{d\mu} (x)} s\xi_\lambda^* (s^{-1} x) - \xi_\lambda^* (x) |^2 ) d\mu(x))^{1/2}\\ &= (4\| \omega\| [\int_X \omega (\frac{ds\mu}{d\mu}(x) s(\xi_\lambda^* (s^{-1} x)) s(\xi_\lambda (s^{-1} x))) d\mu (x)  + \int_X \omega (\xi_\lambda^* (x) \xi_\lambda (x))d\mu(x)\\ & - \int_X \omega (\sqrt{\frac{ds\mu}{d\mu} (x)} s(\xi_\lambda^* (s^{-1} x)) \xi_\lambda (x)) d\mu (x) - \int_X \omega (\xi_\lambda^* (x) \sqrt{\frac{ds\mu}{d\mu} (x) } s( \xi_\lambda (s^{-1} x))) d\mu (x)])^{1/2} \rightarrow 0.
\end{split}
\end{equation*}
The proof is complete since (iii)$\Rightarrow$(iv). $\Box$\\

Suppose that the discrete countable group $G$ acts on the locally compact Hausdorff space $X$ and that $\mu$ is a positive Radon measure on $X$ which is quasi-invariant with respect to this action. We recall from Definition \ref{koop} that the Koopman representation associated to $\mu$ is the group homomorphism $u^\mu : G \rightarrow U(L^2(\mu))$ defined by $\xi \mapsto u_s^\mu (\xi) = \sqrt{ds\mu /d\mu} \xi \circ s^{-1}$ for $\xi \in L^2(\mu)$, $s \in G$. 
\begin{thm}
\label{thm3.5}
Let $G$ be a discrete countable group, $M$ be a commutative von Neumann algebra of the form $M=L^\infty (Y,\nu)$, where $Y$ is a compact Hausdorff space endowed with an action of $G$ and $\nu$ is a finite $G$-quasi-invariant positive Radon meaure on $Y$. Let $X$ be a locally compact Hausdorff space endowed with an action of $G$ and let $\mu$ be a $G$-quasi-invariant positive Radon measure on $X$. Suppose that the Koopman representation associated to $\mu$ is tempered and that there is a norm $1$ equivariant projection $P: L^\infty (X,\mu) \overline{\otimes} M \rightarrow M$. Then there is a norm $1$ equivariant projection $l^\infty (G) \overline{\otimes} M \rightarrow M$.
\end{thm}
\proof Since $\nu$ is finite, we have the continuous embeddings $L^\infty (Y, \nu) \subset L^2(Y,\nu)$ and $L^2 ((X,\mu) , M) \subset L^2((X,\mu), L^2 (Y,\nu))$. Let now $\epsilon >0$, $F \subset G$ finite and $A \subset M_{*,+}$ be given. By Proposition \ref{prop3.4} there is $\xi \in L^2 ((X,\mu),M)_1$ such that 
\begin{equation}
\label{thm3.5.eq1}
\braket{\xi, \xi} =1, \qquad |\omega(\braket{ \sqrt{\frac{ds\mu}{d\mu}} s \xi(s^{-1} \cdot) , \xi} -1 )| < \frac{\epsilon}{3}, \quad \forall \omega \in A, \; s \in F.
\end{equation}
Approximating $\xi$ by a sequence of simple functions $f_k = \sum_{i=1}^{n_k} \chi_{E_i} \otimes m_i$, we find Borel sets $(E_i)_{i \in \mathbb{N}}$ and elements $(m_i)_{j \in \mathbb{N}}$ in $M$ such that the span $B$ of $\{ \sqrt{\frac{ds\mu}{d\mu}}  \chi_{s E_i}  \otimes s m_j\}_{s \in G, i,j \in \mathbb{N}}$ contains all the translates of $\xi$ (under the map $f \mapsto \sqrt{d s\mu /d\mu } \cdot s ( f \circ s^{-1})$, for $s \in G$) in its closure. Let $H$ be the separable $G$-invariant (under the Koopman representation on $L^2(X,\mu)$) subspace of $L^2(X,\mu)$ generated by the span of the set $\{ \sqrt{\frac{ds\mu}{d\mu}} \chi_{s E_i}\}_{s \in G, i \in \mathbb{N}}$. Let also $K$ be the separable closed subspace of $L^2((X,\mu),M)$ generated by elements in $B$. 

Now choose a simple function $\eta=\sum_{i=1}^N \chi_{E_i} \otimes m_i$ in $B$ such that 
\begin{equation*}
\| \xi - \eta\|_{2,\infty} < \delta = \epsilon/( (100 \|\xi\|_{2,\infty} + 100)\max_{\omega \in A} \| \omega\|).
\end{equation*}
Denote by $s \mapsto u_s=u_s^\mu$ the Koopman representation associated to $\mu$; since it is tempered, by Voiculescu's Theorem we can find an isometry $V: H \rightarrow l^2(G)$ such that 
 \begin{equation*}
 \|\lambda_s V h  - V u_s h \|_2 < \epsilon / (3 N^2 \max_{\omega \in A} \| \omega\| \max_{i,j=1,...,N}\{ \|\chi_i\|_2 \|\chi_j\|_2 \|m_i\| \|m_j\|\}) \quad \mbox{ for every }s \in F,  \; h \in H_1. 
 \end{equation*}
If $W  : H \rightarrow l^2(G)$ is an isometry, then a computation on the simple functions generating $K$ shows that $W \otimes 1 : K \rightarrow l^2(G, M)$ is a well defined contractive linear map (for the associated $2-\infty$-norms) satisfying
 \begin{equation}
 \label{thm3.5.eq2}
 \langle (W \otimes 1) \zeta, (W \otimes 1) \zeta' \rangle_M = \langle \zeta, \zeta' \rangle_M \qquad \mbox{ for all } \zeta, \zeta' \in  K.
 \end{equation}
 The linear map $W \otimes 1$ commutes with the action of $G$ on $M$.\\
It follows from \ref{thm3.5.eq1} and \ref{thm3.5.eq2} that the vector $(V \otimes 1) \xi \in l^2(G, M)$ satisfies 
\begin{equation*}
\langle (V \otimes 1) \xi, (V \otimes 1) \xi \rangle_M = 1.
\end{equation*}
Moreover, for every $\omega \in A$ and $s \in F$ we have
\begin{equation*}
\begin{split}
|\omega& (\braket{(\lambda_s  \otimes s)(( V \otimes 1) \xi), (V \otimes 1) \xi} - 1)| \\ &\leq | \omega( \braket{ (1 \otimes s)((\lambda_s V \otimes 1) \xi) - (1 \otimes s)((V u_s \otimes 1) \xi), (V \otimes 1) \xi}| \\ &+ |\omega (\braket{(1 \otimes s)((V u_s \otimes 1) \xi), (V \otimes 1) \xi }- 1)|\\
&=  | \omega( \braket{ (1 \otimes s)((\lambda_s V \otimes 1) \xi) - (1 \otimes s)((V u_s \otimes 1) \xi), (V \otimes 1) \xi}| \\ 
&+ |\omega (\braket{(V u_s \otimes 1) (1 \otimes s) \xi, (V \otimes 1) \xi }- 1)|\\
&=  | \omega( \braket{ (1 \otimes s)((\lambda_s V \otimes 1) \xi) - (1 \otimes s)((V u_s \otimes 1) \xi), (V \otimes 1) \xi}| \\ 
&+ |\omega (\braket{ (1 \otimes s) ((u_s \otimes 1) \xi),  \xi }- 1)| \leq | \omega( \braket{ (1 \otimes s)[\lambda_s V \otimes 1-V u_s \otimes 1] \xi), (V \otimes 1) \xi}|  + \frac{\epsilon}{3}.
\end{split}
\end{equation*}
Now we compute
\begin{equation*}
\begin{split}
\| & \braket{ (1\otimes s)((\lambda_s V \otimes 1-V u_s \otimes 1) \xi), (V \otimes 1) \xi} -  \braket{ (1 \otimes s)((\lambda_s V \otimes 1-V u_s \otimes 1) \eta), (V \otimes 1) \eta}\|\\
& \leq \| \braket{ (1 \otimes s)((\lambda_s V \otimes 1-V u_s \otimes 1) \xi), (V \otimes 1) \xi}  - \braket{ (1 \otimes s)((\lambda_s V \otimes 1-V u_s \otimes 1) \eta), (V \otimes 1) \xi} \| \\
& + \|\braket{ (1 \otimes s)((\lambda_s V \otimes 1-V u_s \otimes 1) \eta), (V \otimes 1) \xi} - \braket{ (1 \otimes s)((\lambda_s V \otimes 1-V u_s \otimes 1) \eta), (V \otimes 1) \eta}\|\\
&= \|\braket{ (1 \otimes s)((\lambda_s V \otimes 1-V u_s \otimes 1) (\xi - \eta), (V \otimes 1) \xi}\|\\ & + \|\braket{ ( 1\otimes s)((\lambda_s V \otimes 1-V u_s \otimes 1) \eta), (V \otimes 1) (\xi -\eta)}\|\\
& \leq \|(1 \otimes s)((\lambda_s V \otimes 1-V u_s \otimes 1) (\xi - \eta)\|_{2,\infty} \|(V \otimes 1) \xi\|_{2,\infty} \\
&+ \|(1 \otimes s)((\lambda_s V \otimes 1-V u_s \otimes 1) \eta\|_{2,\infty} \|(V \otimes 1) (\xi -\eta)\|_{2,\infty}\\
& \leq 2 \| \xi- \eta\|_{2, \infty} \|\xi\|_{2, \infty}+ 2 \|\eta\|_{2, \infty} \|\xi-\eta\|_{2, \infty} <  4 \delta \|\xi\|_{2, \infty} + 2 \delta^2 <  4 \delta \|\xi\|_{2, \infty} + 2 \delta\\ & < \frac{4 \epsilon}{100 \max_{\omega \in A} \| \omega\|} + \frac{2\epsilon}{100 \max_{\omega \in A} \| \omega\|} < \frac{\epsilon}{3\max_{\omega \in A} \| \omega\|}.
\end{split}
\end{equation*}
Hence, for every $\omega \in A$ and $s \in F$ we have
\begin{equation*}
\begin{split}
| \omega( &\braket{ (1 \otimes s)((\lambda_s V \otimes 1-V u_s \otimes 1) \xi), (V \otimes 1) \xi}| \\
&\leq |\omega (\braket{ (1 \otimes s)((\lambda_s V \otimes 1-V u_s \otimes 1) \xi), (V \otimes 1) \xi} - \braket{ (1 \otimes s)((\lambda_s V \otimes 1-V u_s \otimes 1) \eta), (V \otimes 1) \eta}) | \\
&+ |\omega (\braket{ (1 \otimes s)((\lambda_s V \otimes 1-V u_s \otimes 1) \eta),(V \otimes 1) \eta})|\\
& \leq \| \omega\| \frac{\epsilon}{3 \max_{\omega \in A} \| \omega\|} + \| \omega\| \| \braket{ (1 \otimes s)((\lambda_s V \otimes 1-V u_s \otimes 1) \eta), (V \otimes 1) \eta}\|.
\end{split}
\end{equation*}
Now,
\begin{equation*}
\begin{split}
 \| &\braket{ (1 \otimes s)((\lambda_s V \otimes 1-V u_s \otimes 1) \eta), (V \otimes 1) \eta}\| = \| \sum_{i,j=1}^N \braket{ ( \lambda_s V - V u_s) \chi_{E_i} , V \chi_{E_j}}_{l^2(G)} (s m_i)^* m_j\| \\
 & \leq \sum_{i,j= 1}^N | \braket{( \lambda_s V - V u_s) \chi_{E_i} , V \chi_{E_j}}_{l^2(G)}| \|m_i\| \|m_j\| \\
 &\leq \sum_{i,j=1}^N \| \lambda_s V - V u_s\| \| \chi_{E_i}\|_2 \|\chi_{E_j}\|_2 \|m_i \| \|m_j\| < \frac{\epsilon}{3\max_{\omega \in A} \| \omega\|}.
 \end{split}
\end{equation*}
Thus we obtain
\begin{equation*}
|\omega (\braket{(\lambda_s \otimes s)(( V \otimes 1) \xi), (V \otimes 1) \xi} - 1)| < \frac{\epsilon}{3} +  \frac{\epsilon}{3} +  \frac{\epsilon}{3} = \epsilon \qquad \mbox{ for every } \omega \in A \; s \in F.
\end{equation*}
The result follows from Proposition \ref{prop3.4}. $\Box$

The proof of the following result is the same as the one in \cite{AD1} Proposition 4.3.2 and uses the fact that $l^\infty (G, M) = l^\infty (G) \overline{\otimes} M$ if $G$ is discrete countable.

\begin{cor}
\label{cor1.4}
Let $G$ be a discrete countable group, $M$ be a commutative von Neumann algebra of the form $M=L^\infty (Y,\mu)$, where $Y$ is a compact Hausdorff space endowed with an action of $G$ and $\mu$ is a finite $G$-quasi-invariant positive Radon meaure on $Y$. Then $\mu$ is Zimmer-amenable if and only the pair $(M \overline{\otimes} M, M)$ is amenable and the Koopman representation $u^\mu$ is tempered.
\end{cor}
\proof We need to show that if $\mu$ is Zimmer-amenable, then $(M \overline{\otimes} M, M)$ is an amenable pair. Let $E M \overline{\otimes} M \rightarrow M$ be a projection. The assignment $G \times M \overline{\otimes} M \ni (s, \phi) \mapsto \{s \mapsto s (E(s^{1} \phi))\}$ defines a $G$-equivariant projection $M \overline{\otimes} M \rightarrow l^\infty (G) \overline{\otimes} M$. If $P : l^\infty (G) \overline{\otimes} M \rightarrow M$ is an equivariant projection implementing the Zimmer-amenability of $\mu$, the composition $P \circ E$ implements the amenability of the pair $(M \overline{\otimes} M, M)$. $\Box$

\subsection{Equivalence between measurewise amenability and topological amenability (non-separable case, countable group)}
The main result of this section Theorem \ref{thmmt} is well known (see for example \cite{AD4} Theorem 3.3.7). We include a proof for completeness. The next Lemma will be used also in the forthcoming sections. Following \cite{BEW}, given a covariant representation of a $G$-$C^*$-algebra, we will say that it is cyclic if the associated crossed product representation is cyclic. We will make use of the following 
\begin{lem}
\label{lem6.1}
Let $X$ be a compact Hausdorff space endowed with an action of $G$ and $\mu$ a Radon probability measure on $X$. Then $\mu$ is quasi-invariant if and only if the following holds:
\begin{equation*}
\forall g \in G, f_1 \geq f_2 \geq ..., \; f_n: X \rightarrow [0,1] \; \mbox{cont.} \qquad \inf_n \mu(f_n) >0 \Rightarrow \inf_n \mu(g f_n) >0.
\end{equation*}
\end{lem}
\proof Suppose $\mu$ is quasi-invariant and let such a sequence $\{f_n\}$ be given with $m:= \inf_n \mu(f_n) >0$. Suppose there is $g \in G$ such that $\inf_n \mu(g f_n) =0$. We observe that there are a real positive number $c$ and a Borel set $E$ of positive measure such that $f_n|_E \geq c$ for every $n$, indeed suppose this is not the case, then for every $k \in \mathbb{N}$ there is $n \in \mathbb{N}$ such that $\{ f_n^{-1} ([1/k,1])\}$ has measure zero, but then for $k$ large enough so that $1/k <m$ there is $N$ such that we have that $\mu (\{f_n^{-1} [0,1/k]\})=1$ for every $n >N$, contradicting the fact that $m=\inf_n \mu(f_n)$. So let $c>0$ and $E$ be given satisfying  $f_n|_E \geq c$ for every $n$. Then we have $c^{-1} \mu (f_n) \geq \mu(E)$ eventually and so $0= c^{-1} \inf_n (\mu (g f_n)) \geq \mu (gE)$, a contradiction.\\
On the other hand suppose that the above property is satisfied and the measure is not quasi-invariant. Then there are $g \in G$ and $E$ Borel such that $\mu(E)>0$, $\mu(gE) =0$.  Since the measure is Radon and $\mu(gE)=0$, there are open sets $U_n \supset gE$ such that $\mu(U_n) <1/n$ for every $n$. Let now $K \subset E$ be a compact set such that $\mu(K)  >\mu (E) /2$ and let $f_1 \geq f_2 \geq...$, $f_n :X \rightarrow [0,1]$ continuous with $f_n|_{K} =1$, $\supp (f_n)^\circ \subset g^{-1}(U_n)$ for every $n$ and $\inf_n \mu(f_n) = \mu(K)$ (such a sequence exists by regularity of the measure). Then $\inf_n \mu (f_n) >0$, but $\inf_n \mu(g f_n) \leq \inf_n \mu(U_n) =0$, a contradiction. $\Box$

\begin{thm}
\label{thmmt}
	Let $X$ be a compact Hausdorff space (not necessirarly separable) endowed with the action of a countable discrete group $G$. Let $(\pi, u)$ be a cyclic covariant representation of $(C(X), G)$. Then there is a quasi-invariant Radon probability measure $\nu$ on $X$ such that $L^\infty (X,\nu)$ admits a $\sigma$-weakly continuous $G$-equivariant unital homomorphism into $\pi (C(X))''$. 
\end{thm}
\proof Let $\pi \rtimes u$ be the corresponding representation of $C(X) \rtimes G$ on $H$ and $\xi$ be the cyclic vector. The Hilbert subspace $H_\xi= \overline{\pi(C(X)) \xi}$ is the GNS representation of $C(X)$ associated to the Radon measure $f \mapsto \braket{\pi(f) \xi, \xi}$. First we shall see that $H$ is a \underline{countable} direct sum of cyclic representations of $C(X)$. In order to see this, observe that for every $s \in G$, the Hilbert space $u_s H_\xi$ is cyclic for $C(X)$, with associated cyclic vector $u_s \xi$. Since $\xi$ is a cyclic vector for $C(X) \rtimes G$, it follows that $\overline{\spann (\bigcup_{s \in G} u_s H_\xi)} =H$. Fix an enumeration $\{s_n\} $ of $G$ with $s_0 = e$. For every $n \in \mathbb{N}$ consider the orthogonal decomposition $\spann(\bigcup_{i=0}^n u_{s_i} H_\xi) = \oplus_{i=0}^n H_{i}$, where for every $i$, $H_{i}$ is defined inductively by $H_0 = H_\xi$, $H_{i+1} := \spann (\bigcup_{j=0}^{i+1} u_{s_j} H_\xi ) \ominus H_i$. We will see that each $H_{i}$ is a cyclic representation for $C(X)$. For note that for every $n$ the Hilbert space $\spann (\bigcup_{i=0}^n u_{s_i} H_\xi)$ is invariant under $C(X)$, hence the same is true for its orthogonal complement; it follows that the spaces $H_{i}$ are all invariant under $C(X)$. In particular the projections $P_{H_{i}}$ (we assume $0$ is a projection) all belong to the commutant of $\pi (C(X))$ and so $P_{H_{i}} H = P_{H_{i}}  \spann (\bigcup_{j=0}^i u_{s_j} H_\xi) = P_{H_{i}}  u_{s_i} H_\xi = P_{H_{i}} \overline{\pi(C(X)) u_{s_i} \xi} = \overline{ \pi(C(X)) P_{H_{i}} u_{s_i} \xi}$. Thus $P_{H_i} H$ is the GNS Hilbert space for $C(X)$ associated to the Radon measure $\tilde{\mu}_i : f \mapsto \braket{\pi (f) P_{H_i} u_{s_i} \xi, P_{H_i} u_{s_i} \xi}$. Consider now the set $\Lambda$ given by $\{ n \in \mathbb{N} \; | \; P_{H_n} \neq 0\}$. Pick a summable sequence $(\alpha_n)_{n \in \Lambda}$ in $l^1(\mathbb{N}, \mathbb{R})$ of norm one such that $\alpha_n =0$ for every $n \in \mathbb{N} \backslash \Lambda$ and $(\alpha_n)|_\Lambda$ is a sequence of strictly positive numbers, define $\mu := \sum_{n \in \mathbb{N}}\alpha_n  \mu_n$, where for every $n \in \Lambda$ we let $\mu_n := \tilde{\mu}_n /\|\tilde{\mu}_n\|$ and $\mu_n=0$ for $n \notin \Lambda$.\\
\textbf{ Step 1}: $\mu$ is quasi-invariant. We make use of Lemma \ref{lem6.1}. Let $s \in G$, $f_1 \geq f_2 \geq ....$ a sequence of continuous functions on $X$ taking values in $[0,1]$ and suppose that $\inf_k \mu (f_k) =0$ but $\inf_k \mu (f_k \circ s) \neq 0$. Then, up to taking a subsequence, there are $N \in \mathbb{N}$ and $\epsilon >0$ such that $\mu (f_k \circ s) > \epsilon$ for every $k >N$. There is $M \in \mathbb{N}$ such that $\sum_{n >M} \alpha_n \mu_n (f \circ s) < \epsilon /3$ for every continuous $f: X \ra [0,1]$. There is $R \in \mathbb{N}$ such that $u_s P_{H_n} u_{s_n} \xi \in \oplus_{j=0}^R H_j$ for every $n=0,...,M$. Now, for every $0 \leq n \leq M$, $0 \leq i \leq R$ there is $h_{i,n} \in C(X)$ such that $\|P_{H_i} u_s P_{H_n} u_{s_n} \xi - \pi (h_{i,n}) P_{H_i} u_{s_i} \xi \| < \epsilon /(6 (M+1)(R+1))$ and $\| \pi (h_{i,n}) P_{H_i} u_{s_i} \xi\| \leq 1$. Using the covariance condition we compute
\begin{equation*}
\begin{split}
\mu & (f_k \circ s) = \sum_{n=0}^M \alpha_n \braket{ u_{s^{-1}} \pi(f_k) u_s P_{H_n} u_{s_n} \xi, P_{H_n} u_{s_n} \xi} + \sum_{n >M} \alpha_n \mu_n (f_k \circ s)\\
& \leq \sum_{n=0}^M \alpha_n \braket{ u_{s^{-1}} \pi(f_k) u_s P_{H_n} u_{s_n} \xi, P_{H_n} u_{s_n} \xi} + \epsilon /3 \\
&= \sum_{n=0}^M \sum_{i=0}^R \alpha_n \braket{\pi(f_k) P_{H_i}  u_s P_{H_n} u_{s_n} \xi, P_{H_i} u_s P_{H_n} u_{s_n} \xi} + \epsilon /3.
\end{split}
\end{equation*}
Now we estimate
\begin{equation*}
\begin{split}
\sum_{n=0}^M & \sum_{i=0}^R \alpha_n \braket{\pi(f_k) P_{H_i}  u_s P_{H_n} u_{s_n} \xi, P_{H_i} u_s P_{H_n} u_{s_n} \xi}   \\  
& \leq \sum_{n=0}^M \sum_{i=0}^R \alpha_n|\braket{\pi(f_k) (P_{H_i}  u_s P_{H_n} u_{s_n} \xi - \pi(h_{i,n}) P_{H_i}u_{s_i} \xi), P_{H_i} u_s P_{H_n} u_{s_n}\xi}| \\ & +  \sum_{n=0}^M \sum_{i=0}^R \alpha_n|\braket{\pi(f_k) \pi(h_{i,n}) P_{H_i}u_{s_i} \xi, P_{H_i} u_s P_{H_n} u_{s_n} \xi - \pi(h_{i,n}) P_{H_i}u_{s_i} \xi}| \\&+  \max_{i,n} \| h_{i,n}\|^2_\infty \sum_{n=0}^M \sum_{i=0}^R \langle \pi(f_k) P_{H_i}u_{s_i} \xi, P_{H_i}u_{s_i} \xi\rangle \\
& \leq (M+1) (R+1) \cdot \frac{\epsilon}{6 (M+1) (R+1)} + (M+1) (R+1) \frac{\epsilon}{6 (M+1) (R+1)}  \\ & +\max_{i,n} \| h_{i,n}\|^2_\infty \sum_{n=0}^M \sum_{i=0}^R \langle \pi(f_k) P_{H_i}u_{s_i} \xi, P_{H_i}u_{s_i} \xi\rangle < \frac{2\epsilon}{3},\\
\end{split}
\end{equation*}
for $k$ large enough. Hence eventually $\mu  (f_k \circ s)< \epsilon$, which is a contradiction.\\
\textbf{Step 2}: Let $\{e_\lambda\}$ be an orthonormal basis of $l^2(\Lambda)$ and consider the Banach space $L^2 ((X,\mu), l^2 (\Lambda))$. For every $n \in\Lambda$ let $h_n := d\mu_n /d\mu$ be a Borel representative of the Radon-Nikodym derivative of $\mu_n$ with respect to $\mu$. Let $g_n : X \mapsto l^2 (\Lambda)$ be defined as $g_n (x) =  e_n $ for every $x$. This is a Borel function whose class lies in $L^2 ((X,\mu),l^2 (\Lambda))$. Define the map $U: \spann (\bigcup_{n \in \Lambda} u_{s_n} \pi(C(X)) \xi )\rightarrow L^2 ((X,\mu), l^2 (\Lambda))$ given by $U: \sum_{n \in \Lambda} f_n P_{H_n}u_{s_n} \xi \mapsto \sum_{n \in \Lambda} \sqrt{h_n} f_n g_n$. We have
\begin{equation*}
\begin{split}
\langle U &\sum_{n \in \Lambda}  \pi(f_n) u_{s_n} \xi, U \sum_{m \in \Lambda} \pi(l_m) u_{s_m} \xi \rangle_{L^2 ((X,\mu), l^2 (\Lambda))} = \sum_{n,m \in \Lambda} \braket{ f_n g_n, l_m g_m}_{L^2 ((X,\mu), l^2 (\Lambda))} \\ 
&= \sum_{n,m \in \Lambda} \int_X \sqrt{h_n}(x) \sqrt{h_m}(x)f_n (x) l_m (x)  \braket{e_n, e_m}_{l^2(\Lambda)} \\
&= \sum_{n \in \Lambda} h_n (x) f_n (x) l_n (x)  d\mu (x)= \sum_{n \in \Lambda} \int_X f_n (x) l_n (x) d \mu_n (x)\\
&= \sum_{n \in \Lambda} \braket{ \pi(f_n) P_{H_n} u_{s_n} \xi, \pi(l_n) P_{H_n} u_{s_n} \xi}_{H} = \sum_{n,m} \braket{  \pi(f_n) P_{H_n} u_{s_n} \xi, \pi(l_m)P_{H_m} u_{s_m} \xi}_H\\
&= \langle\sum_{n \in \Lambda}  f_n P_{H_n}u_{s_n} \xi, \sum_{m \in \Lambda} \pi(l_m) P_{H_m} u_{s_m} \xi\rangle_H.
\end{split}
\end{equation*}
Hence $U$ is an isometry. Let $\rho : C(X) \rightarrow \mathbb{B}(L^2 ((X,\mu), l^2 (\Lambda)))$ be the representation given by pointwise multiplication; for every $f \in C(X)$ we have $\rho(f)U =U\pi(f)$. Let $P_K : L^2 ((X,\mu), l^2 (\Lambda)) \rightarrow K$ be the projection onto the range of $U$. Since $K$ is identified with a $C(X)$-invariant subspace of $L^2((X,\mu), l^2(\Lambda))$, it follows that $P_K$ commutes with $\rho(C(X))$. The $*$-homomorphism $\pi (C(X)) \rightarrow U \pi(C(X)) U^* = P_K \rho (C(X)) P_K$, $\pi(f) \mapsto U \pi (f) U^*$ extends to a $\sigma$-weakly continuous $*$-isomorphism 
\begin{equation*}
\pi (C(X))'' \rightarrow (P_K \rho (C(X)) P_K)'' = P_K \rho (C(X))'' P_K.
\end{equation*}
 Hence there is a $\sigma$-weakly continuous unital $*$-homomorphism $\phi : \rho (C(X))'' \rightarrow \pi (C(X))''$. The unitary representation of $G$ on $L^2((X,\mu), l^2(\Lambda))$ given by $v_s := u_s^\mu \otimes 1$ implements a covariant representation of $C(X) \rtimes G$. Since the relation $\phi (v_s \rho (f) v_s^*) = u_s \phi (\rho (f)) u_s^*$ holds on the dense subalgebra $\rho (C(X))$, it follows that $\phi$ is equivariant with respect to the adjoint actions $v$ on $\rho (C(X))''$ and $u$ on $\pi (C(X))''$. By amplification there is a unital $G$-equivariant $\sigma$-weakly continuous $*$-homomorphism $L^\infty (X,\mu) \rightarrow \rho(C(X))''$. $\Box$


\begin{cor}
\label{cor1.6}
Let $X$ be a compact Hausdorf space and $G$ a discrete countable group acting on $X$. Then the action of $G$ is topologically amenable if and only if it is measurewise amenable.
\end{cor}
\proof Suppose that the action of $G$ on $X$ is measurewise amenable. By \cite{BEW} Theorem 5.15 and Lemma 3.28, topological amenability is equivalent to the fact that for every cyclic covariant representation $(\pi , u)$ of $C(X) \rtimes G$ the adjoint action of $G$ on $\pi (C(X))''$ is Zimmer-amenable, which in turn is equivalent, in virtue of \cite{BEW} Theorem 3.17 to the existence of a net of compactly supported, norm-continuous, positive type functions $\theta_\lambda : G \rightarrow \pi(C(X))''$ with $\| \theta_\lambda (e)\| \leq 1$ and $\theta_\lambda (g) \rightarrow 1$ $\sigma$-weakly for every $g \in G$. Let $(\pi, u)$ be such a representation; by Theorem \ref{thmmt} there is a quasi-invariant probability measure $\nu$ such that $L^\infty (X,\nu)$ admits a $\sigma$-weakly continuous $G$-equivariant unital homomorphism $\Phi$ into $\pi (C(X))''$. We can precompose $\Phi$ with a net $\theta_\lambda$ given by the Zimmer-amenability of $\nu$ and obtain the Zimmer-amenability of $\pi (C(X))''$. The result follows since topological amenability implies Zimmer-amenability for every Koopman representation, as shown in the proof of \cite{BEW} Theorem 3.27 (1)$\Rightarrow$ (2). $\Box$

\section{Non-standard extensions of Stone-{\v C}ech boundaries}
\label{sec2}

Let $\mathfrak{X}$ be a countable discrete space endowed with an action of a discrete countable group $G$. A \textit{$G$-equivariant compactification of $\mathfrak{X}$} is a compact Hausdorff space $\Delta_\alpha \mathfrak{X}$ endowed with an action of $G$ (by homeomorphisms) such that there is a $G$-equivariant embedding $\mathfrak{X} \rightarrow \Delta_\alpha \mathfrak{X}$ whose image is a dense open set. The \textit{boundary} $\partial_\alpha \mathfrak{X}$ of $\Delta_\alpha \mathfrak{X}$ is by definition $\Delta_\alpha \mathfrak{X} \backslash \mathfrak{X}$, where we identify $\mathfrak{X}$ with its copy inside $\Delta_\alpha \mathfrak{X}$. By Gelfand duality $G$-equivariant compactifications of $\mathfrak{X}$ correspond to $C^*$-subalgebras of $l^\infty (\mathfrak{X})$ containing $c_0 (\mathfrak{X})$. We will denote the Stone-\v{C}ech compactification of $\mathfrak{X}$ by $\Delta_\beta \mathfrak{X}$ and its boundary by $\partial_\beta \mathfrak{X}$, this is the set of free ultrafilters on $\mathfrak{X}$. To every free ultrafilter $\omega \in \partial_\beta \mathbb{N}$ and $G$-equivariant compactification $\Delta_\alpha \mathfrak{X}$ of $\mathfrak{X}$ we can associate the $C^*$-algebra $I_\omega  (\Delta_\alpha \mathfrak{X}):= \{(f_n)_{n \in \mathbb{N}}  \in \prod_{n \in \mathbb{N}} C(\Delta_\alpha \mathfrak{X}) \; | \; \lim_{n \rightarrow \omega} f_n =0\}$, which is an ideal in $\prod_{\mathbb{N}} C(\Delta_\alpha \mathfrak{X})$; the corresponding ultraproduct $C^*$-algebra is $C(\Delta_\alpha \mathfrak{X})_\omega := (\prod_{n \in \mathbb{N}} C(\Delta_\alpha \mathfrak{X})) / I_\omega (\Delta_\alpha \mathfrak{X})$. \\
Let now be given an enumeration $\{x_n\}$ of the discrete countable $G$-space $X$ and an ultrafilter $\omega \in \partial_\beta \mathbb{N}$. For every $n \in \mathbb{N}$ let $c_n (X)$ be the ideal in $C(\Delta_\alpha \mathfrak{X})$ given by functions with support contained in $\{x_0,...,x_n\}$. For every $n \in \mathbb{N}$ let $\mathcal{C}_{n,\alpha} (\mathfrak{X})_\omega:= \prod_{k \in \mathbb{N}} c_n (\mathfrak{X}) /I_\omega (\Delta_\alpha \mathfrak{X})$ be the corresponding ideal in $C(\Delta_\alpha \mathfrak{X})_\omega$. Define the closed $G$-invariant ideal $\mathcal{C}_{0,\alpha}(\mathfrak{X})_\omega := \overline{\bigcup_n \mathcal{C}_{n,\alpha}(\mathfrak{X})}$. The \textit{non-standard extension of $\partial_\alpha \mathfrak{X}$} (relative to $\omega$) is the spectrum of the $C^*$-algebra $C(\Delta_\alpha \mathfrak{X})_\omega /\mathcal{C}_{0,\alpha}(\mathfrak{X})_\omega$ and is denoted $\partial_{\alpha, \omega} \mathfrak{X}$.

Non-standard extensions of boundaries have a simple topological realization.

  
  \begin{lem}
  \label{lem1}
  Let $\omega \in \partial_\beta \mathbb{N}$ and $\psi_\omega : \sigma(C(\Delta_\alpha \mathfrak{X})_\omega) \rightarrow \Delta_\alpha \mathfrak{X}$ be the surjection which is dual to the diagonal embedding. Then $\psi_\omega^{-1} (\partial_\alpha \mathfrak{X}) =\partial_{\alpha,\omega} \mathfrak{X}$.
  \end{lem}
  \proof We need to prove the following: let $A=C(Z) \subset B=C(Y)$ be an inclusion of commutative unital $C^*$-algebras (preserving the identity) and let $\psi: Y \rightarrow Z$ be the associated surjective map, let $U \subset Z$ be a non-empty open set, then $C_0(\psi^{-1}(U))$ equals the ideal $I_U$ generated by $C_0(U)$ inside $C(Y)$. This would suffice since in this case we have $C(\psi^{-1} (Z \backslash U))= C(Y)/C_0 (\psi^{-1}(U)) = C(Y)/I_U$ and in the case at hand, taking $U=X$ and $Y=\sigma (C(\Delta_\alpha X)_\omega)$, we obtain $I_U= \mathcal{C}_{0,\alpha}(X)_\omega$. So let $U$ be such an open set and let $f \in C_0 (\psi^{-1}(U))$; for every $\epsilon >0$ we find a compact set $C \in \psi^{-1} (U)$ such that $\| f|_{Y \backslash C}\| < \epsilon$, let then $g $ be a continuous function on $Z$ with values in $[0,1]$ such that $g|_{\psi (C)} =1$; we have $\| f - f g \circ \psi \| < 2 \epsilon$. Hence every element in $C_0 (\psi^{-1} (U))$ is in the ideal generated by $\psi^* (C_0 (U))$. $\Box$
  
  We specialize now to the case in which the $G$-equivariant compactification of $X$ is the Stone-{\v C}ech compactification.
   
  \begin{lem}
  \label{lem2}
  Let $A \subset B$ be an inclusion of $C^*$-algebras. There is a unique intermediate $C^*$-algebra $A \subset C \subset B$ such that $A$ is an ideal in $C$ and $C$ is maximal among the $C^*$-subalgebras of $B$ containing $A$ as an ideal.
  \end{lem} 
  \proof This is an application of Zorn Lemma. Consider the set $m$ of sub-$C^*$-algebras of $B$ containing $A$ as an ideal with the partial order given by inclusion. Let $c$ be a chain (that is a totally ordered subset of $m$), then $(\bigcup_{C \in c} C) A=A (\bigcup_{C \in c} C)=A$ and so $\overline{(\bigcup_{C \in c} C)}A = A \overline{(\bigcup_{C \in c} C)}=A$. Hence every chain has an upper bound and there is a maximal element $C$. This element has to be unique, for if $C_1$ and $C_2$ are two different such maximal elements, then $C^*(C_1, C_2)$ would be a strict majorant for both $C_1$ and $C_2$. $\Box$
  
For an inclusion $A \subset B$ of $C^*$-algebras we will refer to the unique intermediate $C^*$-algebra $C$ given by the previous Lemma as to the \textit{multiplier algebra of $A$ relative to $B$} and denote it by $M_B(A)$.

  As above we consider a countable discrete space $\mathfrak{X}$ and fix an enumeration $\mathfrak{X}=\{x_n\}_{n \in \mathbb{N}}$. For every $n \in \mathbb{N}$ let $P_n \in \mathbb{B}(l^2(X))$ be the projection on the linear span of $\{\delta_{x_0},..., \delta_{x_n}\}$ and $M_{n+1, \omega} := P_n \mathbb{B}(l^2(\mathfrak{X}))_\omega P_n$ (where $P_n$ is represented diagonally in $\mathbb{B}(l^2(\mathfrak{X}))_\omega$). Define 
  \begin{equation*}
  \mathcal{K}(l^2(\mathfrak{X}))_\omega := \overline{\bigcup_n M_{n,\omega}} \subset \mathbb{B}(l^2(\mathfrak{X}))_\omega \quad \mbox{ and } \quad \mathcal{B}(l^2(\mathfrak{X}))_\omega:= M_{\mathbb{B}(l^2(\mathfrak{X}))_\omega} (\mathcal{K}(l^2(\mathfrak{X}))_\omega).
  \end{equation*}
   Note that the embedding of $l^\infty (\mathfrak{X})_\omega$ inside $\mathbb{B}(l^2(\mathfrak{X}))_\omega$ takes values in $\mathcal{B}(l^2(\mathfrak{X}))_\omega$ and that $\mathcal{K}(l^2(\mathfrak{X}))_\omega \cap l^\infty (\mathfrak{X})_\omega= \mathcal{C}_{0,\beta}(\mathfrak{X})_\omega$.
  
 If $\mathfrak{X}$ is endowed with an action of the discrete countable group $G$, we will denote by $\sigma$ the associated unitary representation on $l^2(\mathfrak{X})$; given $\omega \in \partial_\beta \mathbb{N}$, we denote by $\sigma_\omega$ the unitary representation given by $g \mapsto (\sigma (g))_{n \in \mathbb{N}} + I_\omega\in \mathbb{B}(l^2(\mathfrak{X}))_\omega$. Similarly, we let $\sigma_{\mathbb{N}}$ be the unitary representation of $G$ given by $g \mapsto \sigma_{\mathbb{N}} (g):= (\sigma (g))_{n \in \mathbb{N}} \in \prod_{n \in \mathbb{N}} \mathbb{B}(l^2(\mathfrak{X}))$. Denote by $\mathcal{K}(l^2(\mathfrak{X}))_{\mathbb{N}}$ the $C^*$-algebra given by $\overline{\bigcup_{n \in \mathbb{N}} P_n (\prod_\mathbb{N} \mathbb{B}(l^2(\mathfrak{X}))) P_n} \subset \prod_\mathbb{N} \mathbb{B}(l^2(\mathfrak{X}))$.
 
 Consider the diagonal embedding $i_\omega : l^\infty (\mathfrak{X}) \rightarrow l^\infty (\mathfrak{X})_\omega$. The dual map $\psi_\omega$, when restricted to $\partial_{\beta, \omega} X$, is "nuclear" in a suitable sense. The following Proposition makes this precise and follows the approach used in \cite{Oz} Proposition 4.1.
 
 \begin{prop}
 \label{prop3}
 Let $\mathfrak{X}$ be a discrete countable space endowed with an action of an exact discrete countable group $G$ and let $\omega \in \partial_\beta \mathbb{N}$. The following are equivalent:
 \begin{itemize}
 \item[(i)]: The action of $G$ on $\partial_{\beta,\omega} \mathfrak{X}$ is topologically amenable,
 \item[(ii)]: there is a ucp map $\phi: C^*_\lambda (G) \rightarrow \mathcal{B}(l^2(\mathfrak{X}))_\omega \subset \mathbb{B}(l^2(\mathfrak{X}))_\omega$ such that $\phi(\lambda(g)) - \sigma_\omega(g) \in \mathcal{K}(l^2(\mathfrak{X}))_\omega$ for every $g \in G$,
 \item[(iii)]: there is a nuclear ucp map $\psi: C^*_\lambda (G) \rightarrow \prod_{n \in \mathbb{N}} \mathbb{B}(l^2(\mathfrak{X}))$ such that $\psi (\lambda(g)) - \sigma_{\mathbb{N}} (g) \in \mathcal{K}(l^2(\mathfrak{X}))_\mathbb{N}+I_\omega$ for every $g \in G$,
 \item[(iv)]: there is a nuclear ucp map $\psi: C^*_\lambda (G) \rightarrow \prod_{n \in \mathbb{N}} \mathbb{B}(l^2(\mathfrak{X}))$ such that for every finite set $F \subset G$ and every $\epsilon >0$ there are an $N \in \mathbb{N}$ and an element $A \in \omega$ satisfying  $\|P_N( \psi (\lambda(g)) - \sigma_{\mathbb{N}}(g))_i P_N -( \psi (\lambda(g)) - \sigma_{\mathbb{N}}(g))_i \| < \epsilon$ for every $g \in F$, $i \in A$.
 \item[(v)]: the action on $\partial_\beta \mathfrak{X}$ is topologically amenable.
 \end{itemize}
 \end{prop}
 \proof $(i)\Rightarrow (ii) \Rightarrow (iii)$: In virtue of Lemma \ref{lem1}, $C(\partial_{\beta,\omega} X) = l^\infty(X)_\omega /\mathcal{C}_{0,\beta}(\mathfrak{X})_\omega$ and as we already observed $\mathcal{C}_{0,\beta}(\mathfrak{X})_\omega =  \mathcal{K}(l^2(\mathfrak{X}))_\omega \cap l^\infty (\mathfrak{X})_\omega$. Since $l^\infty (\mathfrak{X})_\omega \subset  \mathcal{B}(l^2(\mathfrak{X}))_\omega$, we conclude that $C(\partial_{\beta,\omega} \mathfrak{X}) \subset \mathcal{B}(l^2(\mathfrak{X}))_\omega / \mathcal{K}(l^2(\mathfrak{X}))_\omega$. For every $g \in G$ the associated unitary $\sigma_\omega (g)$ in $\mathbb{B}(l^2(\mathfrak{X}))_\omega$ lies in $\mathcal{B}(l^2(\mathfrak{X}))_\omega$ and so there is a covariant representation of the pair $(l^\infty(\mathfrak{X})_\omega, G)$ inside $\mathcal{B}(l^2(\mathfrak{X}))_\omega$, which induces a covariant representation of $(C(\partial_{\beta, \omega} \mathfrak{X}), G)$ inside $\mathcal{B}(l^2(\mathfrak{X}))_\omega / \mathcal{K}(l^2(\mathfrak{X}))_\omega$. Since the action is topologically amenable, the $C^*$-algebra $C(\partial_{\beta, \omega} \mathfrak{X})\rtimes_r G$ is nuclear and we have $*$-homomorphisms $C^*_\lambda (G) \rightarrow C(\partial_{\beta, \omega} \mathfrak{X})\rtimes_r G \simeq C(\partial_\beta \mathfrak{X}_\omega)\rtimes G \rightarrow \mathcal{B}(l^2(\mathfrak{X}))_\omega / \mathcal{K}(l^2(\mathfrak{X}))_\omega$. By the Choi-Effros lifting theorem (\cite{CE} Theorem 3.10) we obtain the map $\phi$. Another application of Choi-Effros lifting theorem shows the existence of $\psi$.
 
 $(iii)\Rightarrow (iv)$: this is obtained by unraveling the definitions and using the fact that a finite number of elements of an ultrafilter always have non-empty intersection.
 
 $(iv) \Rightarrow (v)$: Let $F_k \subset F_{k+1} \subset G$ be an exhaustion of $G$ by finite symmetric sets. By (iv) there are sequences $n_k \rightarrow \infty$, $A_k \in \omega$ nested sets such that for every $g \in F_k$ we have $\|P_{n_k} (\psi_i (\lambda (g))-\sigma_{\mathbb{N}}(g))P_{n_k} - (\psi_i (\lambda (g))-\sigma_{\mathbb{N}}(g))\| < 1/k^2$ for every $i \in A_k$ ($\psi_i$ is the projection on the $i$-th factor of $\psi$). By Voiculescu's Theorem (\cite{Da} Theorem II.5.3) for every $k$ there are isometries $V_{k,i} : l^2(\mathfrak{X}) \rightarrow l^2(G)$ for every $i \in A_k$ such that $\psi_i (\lambda (g)) - V_{k,i}^*\lambda(g) V_{k,i} \in \mathbb{K}(l^2(\mathfrak{X}))$ for every $g \in G$. In particular, by choosing for every $k$ any $i \in A_k$, we find a sequence of ucp maps $\psi_k: C^*_\lambda (G) \rightarrow \mathbb{B}(l^2(\mathfrak{X}))$ and of isometries $V_k: l^2(\mathfrak{X}) \rightarrow l^2(G)$ such that $ \psi_k (\lambda (g)) - V_{k}^*\lambda(g) V_{k} \in \mathbb{K}(l^2(\mathfrak{X}))$ for every $k$, $g \in G$ and 
 \begin{equation}
 \label{eq1}
 \|P_{n_k} (\psi_k (\lambda (g))-\sigma(g)) P_{n_k} - (\psi_k (\lambda (g))-\sigma(g))\| < 1/(8k^4) \qquad \mbox{ for every } g \in F_k.
 \end{equation}
  For $x \in X$ let $\mu^{(k)}_x=|V_k \delta_x|^2 \in l^1(G)$. Then for $g \in F_k$ we have
  \begin{equation*}
   \|\mu^{(k)}_{gx} - g\mu^{(k)}_x\|_1 = \sum_{h \in G} | |V_k \delta_{gx}|^2(h)- \lambda(g) |V_k \delta_x|^2 (h)|\leq 2 (2-2\Re (\langle\delta_{gx}, V_k^*\lambda (g) V_k \delta_x\rangle))^{1/2}.
  \end{equation*}
 We want to check that $\braket{\delta_{gx}, V_k^*\lambda (g) V_k \delta_x} = \braket{\sigma(g) \delta_x, V_k^* \lambda (g) V_k}\rightarrow 1$. For every fixed $k$, for every $g \in G$, we have $\|( \psi_k (\lambda(g))  \delta_x - V_k^* \lambda(g)V_k ) \delta_x\| \rightarrow 0$ for $x \rightarrow \infty$ and so there is a finite set $C_k \subset \mathfrak{X}$ such that
 \begin{equation*}
 \begin{split}
\|\sigma (g) \delta_x - V_k^* \lambda (g) V_k \delta_x\|_2 &\leq \|\sigma (g) \delta_x -  (\psi_k (\lambda (g))\delta_x +  (\psi_k (\lambda (g))\delta_x - V_k^* \lambda (g) V_k \delta_x\|_2\\
 &\leq \|\sigma (g) \delta_x -  (\psi_k (\lambda (g))\delta_x\|_2 + \| (\psi_k (\lambda (g))\delta_x - V_k^* \lambda (g) V_k \delta_x\|_2\\ & \overset{(\ref{eq1})}{<}1/(2\sqrt{2} k^4) \qquad \mbox{ for every } g \in F_k, \; x \notin C_k,
 \end{split}
 \end{equation*}
 from which it follows that 
 \begin{equation*}
 |1-\braket{\delta_{gx}, V_k^*\lambda (g) V_k \delta_x}| = |\braket{\delta_{gx}, \delta_{gx} - V_k^*\lambda (g) V_k \delta_x}| < 1/(2\sqrt{2} k^4) \qquad \mbox{ for every } g \in F_k, x \notin C_k.
 \end{equation*}
Up to reindexing, we obtain a sequence $\mu^{(k)} : x \rightarrow \mathcal{P}(G)$ of maps and of finite subsets $C_k \subset C_{k+1} \subset \mathfrak{X}$ such that
\begin{equation*}
\| \mu^{(k)}_{gx}- g\mu^{(k)}_{x} \|_1 < 1/k^2 \qquad \forall g \in F_k, \; x \notin C_k.
\end{equation*}
We can suppose that $F_k C_k \subset C_{k+1}$ for every $k$. Let $D_k := F_k C_k$, then the above estimate is still true after replacing $C_k$ with $D_k$. Let $F$ be a finite subset of $G$, then $F \subset F_{n-1}$ for some $n$. If $x \in D_n \backslash D_{n-1}$, then $x$ cannot be written as $hy$ with $h \in F_{n-1}$, $y \in C_{n-1}$ and so $gx$ does not belong to $C_{n-1}$ for every $g \in F$, in particular $gx$ does not belong to $D_{n-2} = F_{n-2} C_{n-2}$; moreover $FD_n = FF_n C_n \subset FC_{n+1} \subset D_{n+1}$. Hence $F (D_n \backslash D_{n-1}) \subset D_{n+1} \backslash D_{n-2}$. Consider the map $h: \mathfrak{X} \rightarrow \mathbb{N}$, $h(x):= \min \{n >0 \; | \; x \in D_n\}$. On $\mathfrak{X} \backslash D_1$ define
\begin{equation*}
\mu(x)=\frac{1}{h(x)} \sum_{k=1}^{h(x)-1} \mu^{(k)}_x 
\end{equation*}
and extend arbitrarily to $\mathfrak{X}$. We want to show that $\lim_{x \rightarrow \infty} \| g\mu (x) - \mu(gx)\|_1 =0$ for every $g \in G$. Let $g \in G$, then there is $N$ such that $g \in F_{N-2}$ and $g(D_n \backslash D_{n-1}) \subset D_{n+1} \backslash D_{n-2}$ for every $n>N$. Let then $x$ belong to such a $D_n$ and not to $D_{n-1}$. We have the following possibilities depending on where $x$ lands under $g$:
\begin{equation*}
(gx \in D_n) \qquad \|\mu(gx)-g\mu(x)\|_1 = \|\frac{1}{n}\sum_{i=1}^{n-1} \mu^{(i)}_{gx} - \frac{1}{n} \sum_{i=1}^{n-1} g\mu^{(i)}_x\|_1 \leq \sum_k k^{-2} n^{-1}
\end{equation*}
\begin{equation*}
\begin{split}
(gx \in D_{n+1}) \qquad &\|\mu(gx)-g\mu(x)\|_1 = \|\frac{1}{n+1} \sum_{i=1}^n \mu^{(i)}_{gx} - \frac{1}{n} \sum_{i=1}^{n-1} g\mu^{(i)}_x\|_1\\
&=  \|\frac{1}{n+1} \sum_{i=1}^{n} \mu^{(i)}_{gx}-  \frac{1}{n} \sum_{i=1}^{n-1} \mu^{(i)}_{gx} +  \frac{1}{n} \sum_{i=1}^{n-1} \mu^{(i)}_{gx}- \frac{1}{n} \sum_{i=1}^{n-1} g\mu^{(i)}_x\|_1\\
& \leq \frac{1}{n+1}\|\mu^{(n)}_{gx}\|_1 + (n-1) |\frac{1}{n+1} - \frac{1}{n}| +  \| \frac{1}{n} \sum_{i=1}^{n-1} \mu^{(i)}_{gx}- \frac{1}{n} \sum_{i=1}^{n-1} g\mu^{(i)}_x\|_1\\
& \leq \frac{1}{n+1} + \frac{n-1}{n (n+1)}+ \frac{1}{n} \|\sum_{i=1}^{n-1} \mu^{(i)}_{gx}- \sum_{i=1}^{n-1} g\mu^{(i)}_x\|_1\\ &\leq \frac{1}{n+1} + \frac{n-1}{n (n+1)}+ \frac{\sum k^{-2}}{n}
\end{split}
\end{equation*}
\begin{equation*}
\begin{split}
(gx \in D_{n-1}) \qquad &\|\mu(gx)-g\mu(x)\|_1 = \|\frac{1}{n-1} \sum_{i=1}^{n-2} \mu^{(i)}_{gx} - \frac{1}{n} \sum_{i=1}^{n-1} g\mu^{(i)}_x\|_1\\
&= \| \frac{1}{n-1} \sum_{i=1}^{n-2} \mu^{(i)}_{gx} - \frac{1}{n} \sum_{i=1}^{n-2} g\mu^{(i)}_x - \frac{1}{n} g\mu^{(n-1)}_x\|_1\\
&\leq \| \frac{1}{n-1} \sum_{i=1}^{n-2} \mu^{(i)}_{gx} - \frac{1}{n-1} \sum_{i=1}^{n-2} g\mu^{(i)}_x\|_1 \\&+ \|\frac{1}{n-1} \sum_{i=1}^{n-2} g\mu^{(i)}_x -\frac{1}{n} \sum_{i=1}^{n-2} g\mu^{(i)}_x\|_1 + \frac{1}{n}\\
&\leq \frac{1}{n-1} \sum k^{-2} + \frac{1}{n-1} + \frac{1}{n} 
\end{split}
\end{equation*}
Hence $\lim_{x \rightarrow \infty} \| g\mu (x)-\mu (gx)\|_1 =0$ for every $g \in G$. The result follows from \cite{Oz} Proposition 4.1. 

$(v) \Rightarrow (i)$: Follows from the existence of a $G$-equivariant continuous map from $\partial_{\beta, \omega} \mathfrak{X}$ to $\partial_\beta \mathfrak{X}$. $\Box$

\begin{cor}
For every exact discrete countable group $G$, biexactness is equivalent to topological amenability of the action of $G \times G$ on $\partial_\beta G_\omega= \psi_\omega^{-1} (\partial_\beta G_\omega)$ for any $\omega \in \partial_\beta \mathbb{N}$.
\end{cor}

The above result gives a non-standard approach to the boundary amenability of group actions. In the following we will investigate more in details such approach in view of the results contained in the previous section.

\subsection{Uniform perturbations of non-standard points}
\label{sec2.1}
\begin{defn} Let $\mathfrak{X}$ be a discrete countable space endowed with an action of a discrete countable group $G$. Let $G=\{s_i\}_{i \in \mathbb{N}}$ be an enumaration of $G$ with $g_0 =\id$ and $\omega \in \partial_\beta \mathbb{N}$. For every uniformly bounded sequence $\xi=(\xi_n)_{n \in \mathbb{N}}$ in $l^2(\mathfrak{X})$ we denote by $\mu_\xi$ the finite Radon measure on $\sigma (l^\infty (\mathfrak{X})_\omega )$ given by $l^\infty (\mathfrak{X})_\omega \rightarrow \mathbb{C}$, $f=(f_n)^\bullet \mapsto \lim_{n \rightarrow \omega} \langle f_n \xi_n, \xi_n \rangle_{l^2(\mathfrak{X})}$. Define the following subsets of $\mathcal{P}(\sigma(l^\infty (\mathfrak{X})_\omega))$:
\begin{equation*}
\mathrm{NS}_\omega (\mathfrak{X}):= \{ \mu_{(\delta_{x_n})} \; |  \; (x_n)_{n \in \mathbb{N}} \mbox{ a sequence in } \mathfrak{X}\},
\end{equation*}
\begin{equation*}
\widetilde{\mathrm{NS}}_{\omega} (\mathfrak{X}) :=\{ \mu_{(\sum_{i \in \mathbb{N}} \alpha_{i}  (s_i\delta_{(x_n)}))} \in \mathcal{P}(\sigma(l^\infty (\mathfrak{X})_\omega)) \; | \; \sum_{i, j \in \mathbb{N}} \alpha_i \alpha_j < \infty , \; \alpha_i \geq 0 \mbox{ for every } i \in \mathbb{N}\},
\end{equation*}
\begin{equation*}
\begin{split}
\mathrm{NS}(P_f (\mathfrak{X})) :=& \{ \mu \in \mathcal{P}(\sigma(l^\infty (\mathfrak{X})_\omega)) \; | \; \exists M \in \mathbb{N}, \; g \in \prod_{\mathbb{N}}(l^\infty (\mathfrak{X})_+), \; |\supp (g_n)| \leq M\;  \forall n \; : \\ &\; \mu (f) = \lim_{n \rightarrow \omega} \braket{ f_n,g_n}_{l^2 (\mathfrak{X})} \; \forall f=(f_n)^\bullet \in l^\infty (\mathfrak{X})_\omega\}.
\end{split}
\end{equation*}
We will refer to elements of $\mathrm{NS}_\omega (\mathfrak{X})$ as to \textit{non-standard points (on $\mathfrak{X}$)}, to elements of $\widetilde{\mathrm{NS}}_{\omega} (\mathfrak{X})$ as to \textit{uniform approximants of non-standard points (on $\mathfrak{X}$)} and to $\mathrm{NS}(P_f (\mathfrak{X})) $ as to \textit{probability measures with non-standard finite support}. Let $\mathcal{P}^\infty (\mathfrak{X})$ be the set of probability measures on $\mathfrak{X}$ which annihilate $c_0 (\mathfrak{X})$ through its natural embedding in $l^\infty(\mathfrak{X})_\omega$. Define $\mathrm{NS}^\infty_\omega (\mathfrak{X}) = \mathrm{NS}_\omega (\mathfrak{X})\cap \mathcal{P}^\infty (\mathfrak{X})$, $\widetilde{\mathrm{NS}}^\infty_{\omega} (\mathfrak{X})= \widetilde{\mathrm{NS}}_{\omega} (\mathfrak{X}) \cap \mathcal{P}^\infty (\mathfrak{X})$ and $\mathrm{NS}^\infty(P_f (\mathfrak{X})) = \mathrm{NS}(P_f (\mathfrak{X})) \cap \mathcal{P}^\infty (\mathfrak{X})$.
\end{defn}

\begin{oss}
If $(x_n)$ and $(y_n)$ are two sequences of points in $\mathfrak{X}$ and they give rise to the same measure on $l^\infty (\mathfrak{X})_\omega$, then there is $A \in \omega$ such that $x_n = y_n$ for every $n \in A$. For let $f=(f_n)=(\delta_{x_n})$, given $\epsilon >0$ there is $A \in \omega$ such that $| \delta_{y_n} (f_n) - \delta_{x_n} (f_n) | =   |\delta_{y_n} (f_n) - 1| <\epsilon$ for every $n \in A$; since $f_n(y_n) \in \{0,1\}$ for every $n$, it follows that $f(y_n)=1$ for every $n \in A$. In particular the elements of $\mathrm{NS}_\omega (\mathfrak{X})$ can be identified with the points of the non-standard version of $\mathfrak{X}$ associated to $\omega$, ${}^* \mathfrak{X}_\omega$ (cf. \cite{Lind}).
\end{oss}

\begin{oss}For what concerns the definition of $\widetilde{\mathrm{NS}}_{\omega} (\mathfrak{X})$ note that, by the assumption on the sequence $(\alpha_i)_{i \in \mathbb{N}}$, the sequence $(\sum_{i=0}^N \alpha_i (s_i \delta_{x_n}))_{N \in \mathbb{N}}$ in $l^\infty (\mathfrak{X})_\omega$ is uniformly Cauchy in the following sense: for every $\epsilon >0$ there is $N \in \mathbb{N}$ such that $\| \sum_{i=0}^N \alpha_i (s_i \delta_{x_n}) - \sum_{i=0}^M \alpha_i (s_i \delta_{x_n})\|_2 < \epsilon$ for every $M >N$, for every $n \in\mathbb{N}$. The element $(\sum_{i \in \mathbb{N}} \alpha_{i}  (s_i\delta_{(x_n)}))$ is understood to be given by the sequence $\xi=(\xi_n)_{n\in \mathbb{N}}$ defined by $\xi_n = \lim_{N \rightarrow \infty} \sum_{i=0}^N \alpha_{i}  (s_i\delta_{(x_n)})$ in $l^2 (\mathfrak{X})$ for every $n$. The sequence $\xi$ is uniformly bounded in $l^2 (\mathfrak{X})$.
\end{oss}

We observe that a measure of the form $\mu_\xi$ is a probability measure if and only if $\xi$ is represented by a sequence of elements in $l^2 (\mathfrak{X})_{1,+}$ and $\mu \in \mathrm{NS}(P_f (\mathfrak{X})) $ if and only if there is $g= (g_n)$ with uniformly bounded support representing $\mu$ such that $g_n \in l^1 (\mathfrak{X})_{1,+}$ for every $n \in\mathbb{N}$. $\mathrm{NS}(P_f (\mathfrak{X}))$ is convex.

\begin{lem}
\label{lem2.5}
Let $\mathfrak{X}$ be a discrete countable space endowed with an action of a discrete countable group $G$ and $\omega \in \partial_\beta \mathbb{N}$. Then
\begin{equation*}
\overline{\cohull (\mathrm{NS}_\omega (\mathfrak{X}))}^{\| \cdot \|} = \overline{\cohull (\widetilde{\mathrm{NS}}_{\omega}(\mathfrak{X}))}^{\| \cdot \|} = \overline{\mathrm{NS} (P_f (\mathfrak{X}))}^{\| \cdot \|}
\end{equation*}
and 
\begin{equation*}
\overline{\cohull (\mathrm{NS}^\infty_\omega (\mathfrak{X}))}^{\| \cdot \|} = \overline{\cohull (\widetilde{\mathrm{NS}}^\infty_{\omega}(\mathfrak{X}))}^{\| \cdot \|} = \overline{\mathrm{NS}^\infty (P_f (\mathfrak{X}))}^{\| \cdot \|},
\end{equation*}
where the norm is intended to be the natural norm on the dual of $l^\infty (\mathfrak{X})_\omega$.
\end{lem}
\proof It is clear that $\mathrm{NS}^{(\infty)}_\omega (\mathfrak{X}) \subset \widetilde{\mathrm{NS}}^{(\infty)}_{\omega} (\mathfrak{X}) \cap \mathrm{NS}^{(\infty)}(P_f (\mathfrak{X}))$; we will show that $\widetilde{\mathrm{NS}}^{(\infty)}_{\omega} (\mathfrak{X}) \subset \overline{ \mathrm{NS}^{(\infty)}(P_f (\mathfrak{X}))}^{\| \cdot \|} $ and $\mathrm{NS}^{(\infty)}(P_f (\mathfrak{X})) \subset \overline{\cohull (\mathrm{NS}^{(\infty)} (\mathfrak{X}))}^{\| \cdot \|}$. For let $\epsilon >0$, $\mu=\mu_\xi$ with $\xi=(\sum_{i \in \mathbb{N}} \alpha_i (s_i \delta_{x_n})) $ be a probability measure and $\alpha_i \geq 0$ for every $i \in \mathbb{N}$. Since the sequence $\xi$ is uniformly Cauchy in the above sense, there is $N \in \mathbb{N}$ such that $\| \mu - \mu_{\xi_N}\| < \epsilon /2$, where $\xi_N$ is the sequence given by $\xi_{N,n}= \sum_{i=0}^N \alpha_i (s_i \delta_{x_n})$ for every $n \in \mathbb{N}$. Hence $\| \mu - \mu_{\xi_N} / \| \mu_{\xi_N}\|\| < \epsilon$. Since $\mu_{\xi_N}/\| \mu_{\xi_N}\|$ is a probability measure in $\mathrm{NS}(P_f (\mathfrak{X}))$, we obtain the inclusion $\widetilde{\mathrm{NS}}_{\omega} (\mathfrak{X}) \subset \overline{\mathrm{NS}(P_f (\mathfrak{X}))}^{\| \cdot \|}$. If $\mu_\xi$ annihilates $c_0 (\mathfrak{X})$, the same is true for $ \mu_{\xi_N} / \| \mu_{\xi_N}$ and so $\widetilde{\mathrm{NS}}^\infty_{\omega} (\mathfrak{X}) \subset \overline{\mathrm{NS}^\infty(P_f (\mathfrak{X}))}^{\| \cdot \|}$. Let now $\eta=(\eta_n)_{n \in \mathbb{N}}$ be a uniformly bounded sequence in $l^1(\mathfrak{X})_{1,+}$ of elements for which the cardinality of the support has an upper bound $M \in \mathbb{N}$ and let again $\epsilon >0$. Denote by $\mu$ the associated measure. For every $i \in \mathbb{N}$ let $A_i := [i\epsilon /  M , (i+1)\epsilon /  M ) \subset \mathbb{R}$. Let $N$ be smallest natural number such that $N \epsilon /M \geq 1$. The union of the ranges of all the $\eta_n$'s is contained in $\bigsqcup_{i=0}^N A_i$. For every $n \in \mathbb{N}$ and $i=0,..., N$ let $D_{n,i}:= \{ x \in \mathfrak{X} \; | \; \eta_n (x) \in A_i\}$ and let $D_i := \{ n \in \mathbb{N} \; | \; D_{n,i} \neq \emptyset\}$. For every $i =0,..., N$ either $D_i \in \omega$ or $D_i \notin \omega$; hence, taking intersections we can find a subset $F \subset \{ 0,...,N\}$ and $D \in \omega$ such that for every $x \in \mathfrak{X}$ and every $n \in D$ there is $i \in F$ with $\eta_n (x) \in A_i$. Let now $i \in F$ and consider the sets $B_{i,k} := \{ n \in \mathbb{N} \; | \; |\eta_n^{-1} (A_i)|=k\}$ for $k \in \mathbb{N}$. Note that $B_{i,k} =\emptyset$ for every $k > M$, for every $i \in F$ and so for every $i \in F$, $D$ equals the disjoint union of a finite number of sets of the form $B_{i,k}$. Taking the intersection of such sets we find $B \in \omega$ such that $\eta_n  = \sum_{i \in F} \chi_{E_{n,i}} \alpha_{i,n}$ for every $n \in B$, where $E_{n,i} \cap E_{n,j} = \emptyset$ for every $i\neq j$ and for every $n \in D$, $|E_{n,i}| = |E_{m,i}|$ for every $m,n \in B$ and $i \in F$, $\alpha_{i,n} \in A_i$ for every $n \in B$. For every $i \in F$ choose $\alpha_i \in A_i$ and let $\tilde{\eta}_n := \sum_{i \in F} \alpha_i \chi_{E_{n,i}}$. For every $f = (f_n)$ and $n \in B$ we have $| \eta_n (f_n) - \tilde{\eta}_n (f_n)| \leq \|f_n\| M (\epsilon /M) = \|f_n\|\epsilon $ and so $\| \eta_n - \tilde{\eta}_n /\| \tilde{\eta}_n\| \|\leq 2\epsilon$; this shows the inclusion $\mathrm{NS}(P_f (\mathfrak{X})) \subset \overline{\cohull (\mathrm{NS} (\mathfrak{X}))}^{\| \cdot \|}$. For the case in which $\eta$ gives rise to a measure in $\mathrm{NS}^\infty (P_f (\mathfrak{X}))$, we proceed in the following way: write $\tilde{\eta}_n := \sum_{i \in F} \alpha_i \chi_{E_{n,i}} = \sum_{i \in F'} \tilde{\alpha}_i \delta_{(x_n)_i}$ where $F'$ is a finite set and $(x_n)_i$ are sequences of elements of $\mathfrak{X}$. Split $F'$ in two disjoint subsets $F' = H \cup K$, where $H := \{ i \in F' \; | \; \delta_{(x_n)_i} \mbox{ annihilates $c_0 (\mathfrak{X})$}\}$ and $K= F' \backslash H$. Every $\delta_{(x_n)_i}$ with $i \in K$ is eventually concentrated on a finite set and so there is a common finite set $E \subset \mathfrak{X}$ on which every such sequence is eventually concentrated, hence $\|\sum_{i \in K} \tilde{\alpha}_i \delta_{(x_n)_i}\|= \sup_{f \in c_0 (\mathfrak{X})_1} | \sum_{i \in K} \tilde{\alpha}_i \delta_{(x_n)_i} (f)| \leq  \sup_{f \in c_0 (\mathfrak{X})_1}|\sum_{i \in F'} \tilde{\alpha}_i \delta_{(x_n)_i}| \leq \| \sum_{i \in F'} \tilde{\alpha}_i \delta_{(x_n)_i} - \mu \| + \sup_{f \in c_0 (\mathfrak{X})_1} |\mu (f)| \leq 2\epsilon$. Thus $\sum_{i \in H} \tilde{\alpha}_i \delta_{(x_n)_i}$ approximates $\mu$ up to $4\epsilon$. $\Box$

\begin{defn}
\label{def02.2}
We denote by $\overline{F}^{(\infty)}_\omega (\mathfrak{X}) := \overline{\cohull (\mathrm{NS}^{(\infty)}_\omega (\mathfrak{X}))}^{\| \cdot \|} = \overline{\cohull (\widetilde{\mathrm{NS}}^{(\infty)}_{\omega}(\mathfrak{X}))}^{\| \cdot \|} = \overline{\mathrm{NS}^{(\infty)} (P_f (\mathfrak{X}))}^{\| \cdot \|}$ any of the closed set of probability measures appearing in the statement of Lemma \ref{lem2.5} and by $\Q \overline{F}^{(\infty)}_\omega (\mathfrak{X})$ the subset of quasi-invariant probability measures. 
\end{defn}

We observe the following 

\begin{lem}
\label{lem2.6}
Let $A$ be a unital $C^*$-algebra and $I \subset A$ a closed ideal. The set $\mathcal{S}(A,I) :=\{ \phi \in \mathcal{S} (A) \; | \; \phi (I)=0\}$ is closed in the norm topology of $A^*$ and there is a topological affine isomorphism $\mathcal{S}(A,I) \simeq \mathcal{S}(A/I)$.\\ Suppose now that $A$ is a commutative $C^*$-algebra endowed with the action of a countable discrete group $G$ and $I$ is $G$-invariant. Denote by 
\begin{equation*}
\begin{split}
\Q\mathcal{P}(\sigma (A)&, \sigma (I)) :=\\& \{ \phi \in \mathcal{S}(A,I) \; | \; \phi \mbox{ corresponds to a $G$-quasi-invariant probability measure } \}.
\end{split}
\end{equation*}
 Then the above isomorphism restricts to an affine isomorphism of the convex subsets $\Q\mathcal{P} (\sigma (A), \sigma (I)) \simeq \Q \mathcal{P} (\sigma (A/I))$.
\end{lem}
\proof Let $\phi \in \mathcal{S}(A,I)$, then $\phi$ factors through $A/I$ and so $\phi (a) = \phi (a+I)$ for every $a \in A$. Hence for every $y \in I$ we have $\phi (x)=\phi (x+y) \leq \| \phi\| \| x+y\|$ and so $|\phi (x+I)| \leq \inf_{y \in I} \| x+y\| = \|x+I\|$, which proves that the linear functional $a+I \mapsto \phi (a+I)$ has norm less or equal to one; moreover, every positive element in $A/I$ is the image under the quotient map of a positive element in $A$ and $1=\phi (1)= \phi (1+I)$; hence this is a state. The map so obtained $\Phi: \mathcal{S}(A,I) \rightarrow \mathcal{S}(A/I)$ respects convex combinations. $\mathcal{S}(A,I)$ is closed: let $\phi_\lambda \rightarrow \phi$ with $\phi_\lambda \in \mathcal{S}(A,I)$ and suppose $\phi \notin \mathcal{S}(A,I)$, then there is $y \in I$ with $\phi (y) \neq 0$, but $\phi (y) = \lim_\lambda \phi_\lambda (y)=0$. $\Phi$ is injective: suppose $\Phi (\phi) = \Phi (\psi)$, then $\phi (a) = \phi (a+I) = \psi (a+I)= \psi (a)$ for every $a \in A$. $\Phi$ is surjective: let $\phi \in \mathcal{S}(A/I)$ and let $q :  \rightarrow A/I$ be the quotient map. Then $\phi \circ q$ is a positive linear functional satisfying $|\phi \circ q (x)| = |\phi (x+I)| \leq \| x+I\| \leq \| x\|$,  hence $\phi \circ q$ has norm less or equal than one. Now, for every $\epsilon >0$ we find $a+I$ with $\inf_{y \in I}\| a + y\| \leq 1$ and $|\phi (a+I)| > 1-\epsilon$; hence we find $y \in I$ with $\| a+y\| < 1+\epsilon$ and $|\phi \circ q (a+y)| > 1-\epsilon$; hence $z=a+y /\|a+y\|$ satisfies $\|z\|=1$ and $|\phi \circ q(z)|> (1-\epsilon) / \|a+y\| > (1-\epsilon)/(1+\epsilon) > 1-2\epsilon $. \\
For the second part we need to check that the state $\phi \in \mathcal{S}(A,I)$ is quasi-invariant if and only if the same is true for $\Phi (\phi)$. It follows from Lemma \ref{lem6.1} that if $\Phi (\phi)$ is quasi-invariant, then the same is true for $\phi$. Suppose now that $\phi$ is quasi-invariant and that $\Phi (\phi)$ is not. Then there are a decreasing sequence $1 \geq f_1 \geq f_2 \geq ...$ of continuous functions on $\sigma (A/I)$ and $g \in G$ such that $\inf_n \phi (f_n) >0$ and $\inf_n \phi (f_n \circ g)=0$. By Gelfand duality the quotient map is dual to an equivariant embedding $\sigma (A/I) \rightarrow \sigma (A)$. By Tietze's extension theorem we can extend $f_1$ to a positive continuous function $h_1$ on $\sigma (A)$ such that $\|h_1\| \leq \|f_1\| \leq 1$ and $\phi (f_1)=\phi (h_1)$, $\phi (f_1 \circ g) = \phi (h_1 \circ g)$. Suppose now we are given $n \in \mathbb{N}$ and a positive continuous extension $h_n$ of $f_n$ such that $\| h_n\| \leq \|f_n\|$, $\phi (h_n) = \phi (f_n)$ and $\phi (h_n \circ g)= \phi (f_n \circ g)$; let $\tilde{h}_{n+1}$ be a positive continuous extension of $f_{n+1}$ given by Tietze's extension theorem and let $h_{n+1} := \tilde{h}_{n+1} \wedge h_n$. Then $h_{n+1} \leq h_n$, $h_{n+1}|_{\sigma (A/I)} = f_{n+1}$ and $h_{n+1} \circ g|_{\sigma (A/I)} = f_{n+1} \circ g|_{\sigma (A/I)}$, whence $\phi (h_{n+1}) = \phi (f_{n+1})$ and $\phi (h_{n+1} \circ g) = \phi (f_{n+1} \circ g)$. The sequence $(h_n)$ so obtained satisfies $1 \geq h_1 \geq h_2\geq ...$, $\inf_n \phi (h_n) = \inf_n \phi (f_n) >0$ and $\inf_n \phi (h_n \circ g) = \inf_n \phi (f_n \circ g)=0$, contradicting quasi-invariance of $\phi$. $\Box$

Since a measure on $l^\infty (\mathfrak{X})_\omega$ annihilates $c_0 (\mathfrak{X})$ if and only if it annihilates $\mathcal{C}_{0,\beta} (\mathfrak{X})_\omega$, it follows from Lemma \ref{lem2.6} that the sets (together with their convex hulls and closures) $\mathrm{NS}^\infty_\omega (\mathfrak{X})$, $\widetilde{\mathrm{NS}}^\infty_{\omega}(\mathfrak{X}))$ and $\mathrm{NS}^\infty (P_f (\mathfrak{X}))$ can be identified with sets of probability measures on $\partial_{\beta, \omega} (\mathfrak{X})$. In the following we want to study in more detail the regularity properties of the quasi-invariant measures belonging to these sets.

\subsection{Automatic pair amenability in $\cohull (\widetilde{\mathrm{NS}}_{\omega}(\mathfrak{X}))$}
In this section we let $\mathfrak{X}$ be a countable set endowed with the action of a countable discrete group $G$. We also fix a free ultrafilter $\omega \in \partial_\beta \mathbb{N}$. We keep the notation described right after the Introduction.  We will exhibit a specific class of quasi-invariant probability measures on the non-standard extensions of $\partial_\beta \mathfrak{X}$ for which pair amenability is automatic. These are convex combinations of certain uniform approximants of non-standard points.
\begin{lem}
\label{lemza1}
Let $\xi=(\xi_n)_{n \in \mathbb{N}}$ be a sequence in $l^2(\mathfrak{X})_1$ and $\mu_\xi$ the associated Radon probability measure on $\sigma(l^\infty (\mathfrak{X})_\omega)$. Let also $\mu_{\xi \otimes \xi}$ be the measure on $\sigma (l^\infty(\mathfrak{X} \times \mathfrak{X})_\omega)$ associated to $\xi \otimes \xi= (\xi_n \otimes \xi_n)_{n \in \mathbb{N}}$. There is an isometry $U_\xi : L^2(\mu_\xi) \otimes L^2(\mu_\xi) \rightarrow L^2(\mu_{\xi \otimes \xi})$. Under this isometry we have $U_\xi \pi_{\mu_\xi \otimes \mu_\xi} ((f_n) \otimes (g_n)) U_\xi^*= \pi_{\mu_{\xi \otimes \xi}} ((f_n \otimes g_n)) U_\xi U_\xi^*$ for every $f,g \in l^\infty(\mathfrak{X})_\omega$.
\end{lem}
\proof There are embeddings $l^\infty (\mathfrak{X})_\omega \otimes 1 \rightarrow l^\infty (\mathfrak{X} \times \mathfrak{X})_\omega$, $(f_n)\otimes 1 \mapsto (f_n \otimes 1)$ and $1 \otimes l^\infty(\mathfrak{X})_\omega \rightarrow l^\infty (\mathfrak{X} \times \mathfrak{X})_\omega$, $1 \otimes (g_n) \mapsto (1 \otimes g_n)$ with commuting range. By the universal property of the tensor product there is a $*$-homomorphism $\pi: l^\infty (\mathfrak{X})_\omega \otimes l^\infty (\mathfrak{X})_\omega \rightarrow l^\infty (\mathfrak{X} \times \mathfrak{X})_\omega$, sending $(f_n) \otimes (g_n)$ to $(f_n \otimes g_n)$. We have $\mu_{\xi \otimes \xi} \circ \pi = \mu_\xi \otimes \mu_\xi : l^\infty (\mathfrak{X})_\omega \otimes l^\infty (\mathfrak{X})_\omega \rightarrow \mathbb{C} $ and so the GNS-Hilbert space associated to $\mu_\xi \otimes \mu_\xi$ can be identified with a closed subspace of $L^2(\mu_{\xi \otimes \xi})$ via the isometry $U_\xi$. The result follows. $\Box$

\begin{lem}
\label{lemza20}
Let $\mu=\mu_\xi \in \widetilde{NS}_\omega (\mathfrak{X})$ be a probability measure associated to a sequence $\xi=(\xi_n)=(\sum_i \alpha_i s_i \delta_{x_n})$ in $l^2(\mathfrak{X})_{1,+}$ for some sequences $(x_n)$ in $\mathfrak{X}$ and $(\alpha_i)$ in $\mathbb{R}_+$ such that $\alpha_i >0$ for every $i$. Then both $\mu$ and $\mu_{\xi^2}$ are quasi-invariant.
\end{lem}
\proof We use the characterization of quasi-invariance given in Lemma \ref{lem6.1}. Let then $(f^{(j)})$ be a decreasing sequence of positive elements of $l^\infty(\mathfrak{X})_\omega$ such that $\inf_j \mu_{\xi}(f^{(j)}) =0$. We can suppose $f^{(j)}=(f^{(j)}_n)^\bullet$ with $f^{(j)}_n \geq 0$ for every $n$. Using the fact that for every $n$ the sum involved only contains positive elements, we obtain
\begin{equation*}
\begin{split}
\mu_{\xi} (f^{(j)})&= \lim_{n \rightarrow \omega} \braket{f^{(j)}_n, \sum_{i,k} \alpha_i \alpha_k \delta_{g_i x_n} \delta_{g_k x_n}}= \lim_{n\rightarrow \omega} \sum_{x \in \mathfrak{X}} f^{(j)}_n (x) \sum_{i,k} \alpha_i \alpha_k \delta_{g_i x_n}(x) \delta_{g_k x_n}(x)\\
& = \lim_{n \rightarrow \omega} \sum_{i,k} \sum_{x \in \mathfrak{X}} \alpha_i \alpha_k f^{(j)}_n (x) \delta_{g_i x_n}(x) \delta_{g_k x_n}(x) = \lim_{n \rightarrow \omega} \sum_{i,k} \alpha_i \alpha_k \braket{f^{(j)}_n, \delta_{g_i x_n} \delta_{g_k x_n}}\\
& \geq \alpha_i \alpha_k \lim_{n \rightarrow \omega} \braket{f^{(j)}_n, \delta_{g_i x_n} \delta_{g_k x_n}} \qquad \forall i,k.
\end{split}
\end{equation*}
Hence, if $\inf_j \mu_{\xi}(f^{(j)})=0$, then for every $i,k \in \mathbb{N}$ we have
\begin{equation*}
\inf_j \lim_{n \rightarrow \omega} \braket{f^{(j)}_n, \delta_{g_i x_n} \delta_{g_k x_n}}_{l^2(\mathfrak{X})}=0 
\end{equation*}
and so
\begin{equation*}
\inf_j \lim_{n \rightarrow \omega} \braket{f^{(j)}_n \circ g , \delta_{g_i x_n} \delta_{g_k x_n}}_{l^2(\mathfrak{X})}=0  \qquad \forall g \in G, \; i,k
\end{equation*}
As above, for every $g \in G$ we have
\begin{equation*}
\begin{split}
\mu_{\xi}(f^{(j)} \circ g)&= \lim_{n \rightarrow \omega} \braket{f^{(j)}_n  \circ g, \sum_{i,k} \alpha_i \alpha_k \delta_{g_i x_n} \delta_{g_k x_n}}= \lim_{n\rightarrow \omega} \sum_{x \in \mathfrak{X}} f^{(j)}_n \circ g (x) \sum_{i,k} \alpha_i \alpha_k \delta_{g_i x_n}(x) \delta_{g_k x_n}(x)\\
& = \lim_{n \rightarrow \omega} \sum_{i,k} \sum_{x \in \mathfrak{X}} \alpha_i \alpha_k f^{(j)}_n \circ g(x) \delta_{g_i x_n}(x) \delta_{g_k x_n}(x) = \lim_{n \rightarrow \omega} \sum_{i,k} \alpha_i \alpha_k \braket{f^{(j)}_n \circ g, \delta_{g_i x_n} \delta_{g_k x_n}}.
\end{split}
\end{equation*}
Suppose that there is $g  \in G$ such that $\inf_j \mu_{\xi} (f^{(j)} \circ g) =c >0$. Let $M \in \mathbb{N}$ be such that $\sum_{\min\{i,k\}>M} \alpha_i \alpha_k < c/3$ and let $l \in \mathbb{N}$ be such that $\lim_{n \rightarrow \omega} \braket{f^{(l)}_n \circ g , \delta_{g_i x_n} \delta_{g_k x_n}}_{l^2(\mathfrak{X})} < c/(3M^2)$ for every $ i,k \leq M$; there is $A \in \omega$ such that $\braket{f^{(l)}_n \circ g , \delta_{g_i x_n} \delta_{g_k x_n}}_{l^2(\mathfrak{X})} < c/(3M^2)$ for every $n \in A$, $ i,k \leq M$ and since $\inf_j \mu_{\xi} (f^{(j)} \circ g) =c$, there is $B \in \omega$ such that $\sum_{i,k} \alpha_i \alpha_k \braket{f^{(l)}_n \circ g, \delta_{g_i x_n} \delta_{g_k x_n}} > 2c/3$ for every $n \in B$. hence for every $n \in A \cap B$ we have
\begin{equation*}
2c/3 < \sum_{i,k} \alpha_i \alpha_k \braket{f^{(l)}_n \circ g, \delta_{g_i x_n} \delta_{g_k x_n}} \leq  \sum_{i,k \leq M} \alpha_i \alpha_k  \braket{f^{(l)}_n \circ g, \delta_{g_i x_n} \delta_{g_k x_n}} + c/3 < 2c/3,
\end{equation*}
a contradiction. This shows the result for $\mu_\xi$. A similar argument shows the quasi-invariance of $\mu_{\xi^2}$. $\Box$

\begin{lem}
\label{lemza2}
Let $\mu=\sum_{m=1}^l \lambda_m \mu_{\xi_m} \in \cohull (\widetilde{\mathrm{NS}}_{\omega}(\mathfrak{X}))$ be a probability measure such that for every $m$ the sequence $(\alpha^{(m)}_i)$ satisfies $\alpha^{(m)}_i >0$ for every $i$. There are a uniformly bounded sequence $\zeta= (\zeta_n)$ in $l^2(\mathfrak{X})_+$ such that $\mu =\mu_\zeta$ and a unitary isomorphism $L^2(\mu_\zeta)\simeq L^2(\mu_{\zeta^2})$; the transpose map is a unital equivariant $\sigma$-weakly continuous isomorphism $L^\infty (\mu_{\zeta^2}) \simeq L^\infty (\mu_\zeta)$.
\end{lem}
\proof Let $f=(f_n)^\bullet \in l^\infty(\mathfrak{X})_\omega$ and $\delta >0$. For every $m=1,..., l$ let $\xi_m = (\sum_i \alpha^{(m)}_{i} \delta_{g_i x^{(m)}_n})$. Then $\mu =\mu_\zeta$ for $\zeta=\sqrt{\sum_{m=1}^l \lambda_m \xi^2_m}$. Let $N \in \mathbb{N}$ be such that $\sum_{i,k >N} \alpha^{(m)}_{i} \alpha^{(m)}_{k} < \delta/ \|f\|^2_\infty$ for every $1 \leq m \leq l$. The sequence $(f_n \chi_{\bigcup_{m=1}^l\{g_i x^{(m)}_n\}_{i=0}^N} / \zeta_n)$ is uniformly bounded. Moreover,
\begin{equation*}
\begin{split}
\|f &- (\frac{f_n \chi_{\bigcup_{m=1}^l\{g_i x^{(m)}_n\}_{i=0}^N}}{\zeta_n } \zeta_n)\|^2_{2, \mu_\zeta} = \lim_{n \rightarrow \omega} \sum_{x \notin \bigcup_{m=1}^l\{g_i x^{(m)}_n\}_{i=0}^N} |f_n(x)|^2 (\sum_{m=1}^l \lambda_m \xi_m^2 (x))\\& \leq \|f\|^2_\infty \sum_{m=1}^l \lambda_m \sum_{i,j>N} \alpha^{(m)}_{i} \alpha^{(m)}_{j} <  \delta.
\end{split}
\end{equation*}
Hence $\eta_{\mu_\zeta} (l^\infty (\mathfrak{X})_\omega \cdot \zeta)$ is a dense linear subspace of $L^2 (\mu_\zeta)$. As a consequence the isometry $U: L^2(\mu_{\zeta^2}) \rightarrow L^2(\mu_\zeta)$, which sends $\eta_{\mu_{\zeta^2}}(f) \in \eta_{\mu_{\zeta^2}}(l^\infty (\mathfrak{X})_\omega)$ to $ \eta_{\mu_\zeta}(f \zeta)$ is a unitary. 
For every $f \in l^\infty (\mathfrak{X})_\omega$ represented as a multiplication operator on $L^2(\mu_{\zeta^2})$, we have $U\pi_{\mu_{\zeta^2}}(f) = \pi_{\mu_{\zeta} }(f) U$ and so $U\pi_{\mu_{\zeta^2}}(f) U^* = \pi_{\mu_{\zeta} }(f)$. Hence the map $T \mapsto UTU^*$ implements a unital equivariant $\sigma$-weakly continuous isomorphism $L^\infty (\mu_{\zeta^2}) \simeq L^\infty (\mu_\zeta)$. $\Box$

\begin{lem}
\label{lemza3}
Let $\xi=(\xi_n)_{n \in \mathbb{N}}$ be a sequence in $l^2(\mathfrak{X})_{1, +}$ representing a measure $\mu_\xi$ which is a convex combination of uniform approximants of non-standard points $\mu_{(\sum_i \alpha_i^{(m)} \delta_{g_i x_n^{(m)}})}$ associated to sequences $(\alpha^{(m)}_i)$ satisfying $\alpha^{(m)}_i>0$ for every $i$, $m$. Then $\eta_{\mu_{\xi \otimes \xi}} (f)$ belongs to the $L^2(\mu_{\xi \otimes \xi})$-closure of $\eta_{\mu_{\xi \otimes \xi}} (\pi (l^\infty (\mathfrak{X})_\omega \otimes l^\infty (\mathfrak{X})_\omega))$ for every $f \in l^\infty (\mathfrak{X} \times \mathfrak{X})_\omega$, where $\pi$ is as in the proof of Lemma \ref{lemza1}.
\end{lem}
\proof We first observe that for every $\epsilon >0$ there is a finite cardinality $M$ such that for every $n \in \mathbb{N}$ there is a set $E_n \subset \mathfrak{X} \times \mathfrak{X}$ with $|E_n|=M$ satisfying $\sum_{(x,y) \in \mathfrak{X} \times \mathfrak{X}  \backslash E_n} \xi_n^2 (x) \xi_n^2 (y) < \epsilon$. To this end it is enough to show that given such $\epsilon$ there is a finite cardinality $N$ such that for every $n \in \mathbb{N}$ there is a subset $F_n \subset \mathfrak{X}$ with $|F_n| = N$ satisfying $\sum_{x \in X \backslash F_n} \xi_n^2(x) < \epsilon /2$ (hence $E_n = F_n \times F_n$ will work). To this end, as in the proof of Lemma \ref{lemza2}, we can take $N' \in \mathbb{N}$ such that $\sum_{i,k >N'} \alpha^{(m)}_{i} \alpha^{(m)}_{k} < \epsilon /2$ and the sets $F_n := \bigcup_{l=1}^m \{ g_i x_n^{(m)}\}_{i=0}^{N'}$ satisfy $\sum_{x \in X \backslash F_n} \xi_n^2(x) < \epsilon /2$ for every $n$. Hence for every $f=(f_n)^\bullet \in l^\infty(\mathfrak{X} \times \mathfrak{X})_\omega$ we can find a sequence of sets $E_n \subset \mathfrak{X} \times \mathfrak{X}$ of uniformly bounded cardinality $M$ such that  $\| f - (f_n \chi_{E_n})\|^2_{2, \mu_{\xi \otimes \xi}} < \epsilon$. The result follows since for every $n \in \mathbb{N}$ we have $f_n \chi_{E_n} = \sum_{(x,y) \in E_n} f_n (x,y) \delta_x \otimes \delta_y$. $\Box$\\

\begin{oss}
\label{remza4}
It follows from Lemma \ref{lemza3} that in the case of measures which are convex combinations of uniform approximants of non-standard points the isometry of Lemma \ref{lemza1} is a unitary.\\
\end{oss}

\begin{prop}
\label{propza5}
Let $\mu$ be a quasi-invariant Radon probability measure on $\sigma(l^\infty (\mathfrak{X})_\omega)$ which is a convex combination of uniform approximants of non-standard points associated to sequences $(\alpha^{(m)}_i)$ satisfying $\alpha^{(m)}_i>0$ for every $i$, $m$. The pair $(\sigma(l^\infty (\mathfrak{X})_\omega) \times \sigma(l^\infty (\mathfrak{X})_\omega), \mu \otimes \mu), (\sigma(l^\infty (\mathfrak{X})_\omega) , \mu))$ is amenable.
\end{prop}
\proof Let $\xi=(\xi_n)$ be the sequence of vectors in $l^2(\mathfrak{X})_+$ implementing $\mu= \mu_\xi$ (cf. the proof of Lemma \ref{lemza1}). By Lemma \ref{lemza1} and Remark \ref{remza4} there is a unitary 
\begin{equation*}
U_\xi : L^2(\mu_\xi) \otimes L^2(\mu_\xi) \rightarrow L^2(\mu_{\xi \otimes \xi}) \qquad \mbox{ such that }
\end{equation*}
\begin{equation*}
U_\xi \pi_{\mu_\xi \otimes \mu_\xi} ((f_n) \otimes (g_n)) U_\xi^*= \pi_{\mu_{\xi \otimes \xi}} ((f_n \otimes g_n))  \qquad \forall f=(f_n)^\bullet, \; g=(g_n)^\bullet \in l^\infty (\mathfrak{X})_\omega.
\end{equation*}
By Lemma \ref{lemza2} there is a unitary
\begin{equation*}
V_\xi : L^2(\mu_{\xi^2}) \rightarrow L^2(\mu_\xi) \qquad \mbox{ such that }
\end{equation*}
\begin{equation*}
V_\xi \pi_{\mu_{\xi^2}}(f) V_\xi^* = \pi_{\mu_{\xi} }(f) \qquad \forall f=(f_n)^\bullet \in l^\infty (\mathfrak{X})_\omega.
\end{equation*}
Consider now the "diagonal" $*$-homomorphism $\hat{\pi} : l^\infty (\mathfrak{X})_\omega \rightarrow l^\infty (\mathfrak{X} \times \mathfrak{X})_\omega$ given by $\hat{\pi} ((f_n)^\bullet) = (g_n)^\bullet$, with $g_n (x,y) = \delta_{x,y} f_n (x)$. This induces an isometry 
\begin{equation*}
W_\xi : L^2(\mu_{\xi^2}) \rightarrow L^2(\mu_{\xi \otimes \xi})
\end{equation*}
Let now $\check{\pi} : l^\infty (\mathfrak{X} \times \mathfrak{X})_\omega \rightarrow l^\infty (\mathfrak{X})_\omega$ be the "restriction" $*$-homomorphism given by $(\check{\pi} (f))_n= (h_n)^\bullet$ with $h_n (x) = f_n (x,x)$ for $f=(f_n)^\bullet$. We have
\begin{equation*}
W^*_\xi \pi_{\mu_{\xi \otimes \xi}} (f) W_\xi = \pi_{\mu_{\xi^2}} (\check{\pi}(f))  \qquad \forall f =(f_n)^\bullet \in l^\infty (\mathfrak{X} \times \mathfrak{X})_\omega.
\end{equation*}
Consider now the $\sigma$-weakly continuous unital linear map 
\begin{equation*}
\gamma : \mathbb{B}(L^2 (\mu_\xi) \otimes L^2(\mu_\xi)) \rightarrow \mathbb{B}(L^2(\mu_\xi))\qquad \mbox{ given by }
\end{equation*}
\begin{equation*}
a \mapsto V_\xi W^*_\xi U_\xi a U_\xi^* W_\xi V_\xi^*
\end{equation*}
 and let $f= (f_n)^\bullet, g=(g_n)^\bullet \in l^\infty (\mathfrak{X})_\omega$, $v \in L^2(\mu_\xi)$. We obtain
\begin{equation*}
\begin{split}
V_\xi W^*_\xi U_\xi &\pi_{\mu_\xi \otimes \mu_\xi} (f\otimes g) U_\xi^* W_\xi V_\xi^* v = V_\xi W^*_\xi \pi_{\mu_{\xi \otimes \xi}} (\pi (f \otimes g))  W_\xi V_\xi^* v\\
&=V_\xi W^*_\xi \pi_{\mu_{\xi \otimes \xi}} (\pi (f \otimes g)) W_\xi V_\xi^* v\\
&= V_\xi \pi_{\mu_{\xi^2}} (\check{\pi}(\pi (f \otimes g)))V_\xi^* v = \pi_{\mu_\xi} (\check{\pi} \circ \pi (f \otimes g)) v.
\end{split}
\end{equation*}
The $*$-homomorphism $\check{\pi} \circ \pi : l^\infty (X)_\omega \otimes l^\infty (X)_\omega \rightarrow l^\infty (X)_\omega$ is surjective and equivariant (on $l^\infty (X)_\omega \otimes l^\infty (X)_\omega$ we consider the diagonal action of $G$). Hence $\gamma$ gives an equivariant ($\sigma$-weakly continuous) projection  $L^\infty (\mu_\xi \otimes \mu_\xi) \rightarrow L^\infty (\mu_\xi)$ which satisfies $\gamma (f \otimes 1)= f$ for every $f \in L^\infty (\mu_\xi)$. $\Box$

\subsection{Internal measures and Calkin states}
\label{susec2.3}
In this subsetion we fix again a countable set $\mathfrak{X}$ endowed with the action of a discrete countable group $G = \{ s_i\}$ with $s_0=e$. We will not fix a free ultrafilter $\omega \in \partial_\beta \mathbb{N}$. We will show that regularity properties of the Calkin representation $C^*(G) \rightarrow \mathbb{B}(l^2(\mathfrak{X}))/\mathbb{K}(l^2(\mathfrak{X}))$ are encoded in regularity properties of the Koopman representations associated to a specific class of quasi-invariant probability measures on non-standard extensions of the Stone-{\v C}ech boundary of $\mathfrak{X}$.

\begin{defn}
For every $\omega \in \partial_\beta \mathbb{N}$ let $\mathcal{P}(\mathfrak{X})_\omega$ (resp. $\mathcal{P}^\infty (\mathfrak{X})_\omega$) be the set of probability measures on $\sigma(l^\infty (\mathfrak{X})_\omega)$ (resp. $\partial_{\beta, \omega} (\mathfrak{X})$) given by sequences of probability measures $(\mu_n)$ on $\mathfrak{X}$ via the formula $f \mapsto \lim_{n \rightarrow \omega} \mu_n (f_n)$ for $f=(f_n)^\bullet$. Denote by $\Q \mathcal{P}(\mathfrak{X})_\omega$ (resp. $\Q \mathcal{P}^\infty(\mathfrak{X})_\omega$) the subset of $G$-quasi-invariant probability measures in $\mathcal{P}(\mathfrak{X})_\omega$ (resp. $\mathcal{P}^\infty(\mathfrak{X})_\omega$). We will refer to elements in $\mathcal{P}(\mathfrak{X})_\omega$ as to \textit{internal probability measures}.
\end{defn}
Every measure in $\mathcal{P}(\mathfrak{X})_\omega$ can be written as a measure of the form $\mu_\xi$ for some sequence $\xi =(\xi_n)$ of elements in $l^2(\mathfrak{X})_{1,+}$.

\begin{defn}
\label{def2.2}
Let $\omega \in \partial_\beta \mathbb{N}$, $\tilde{\xi}=(\tilde{\xi}_n)$ be a sequence in $l^2(\mathfrak{X})_1$, $\tilde{\mu}=\mu_{\tilde{\xi}}$ the associated probability measure and $(\alpha_i)$ be a sequence of strictly positive real numbers such that $\sum_{i,j} \alpha_i \alpha_j < \infty$. Define the probability measure on $l^\infty (\mathfrak{X})_\omega$ by $\mu=\mu_{\xi}$, where $\xi=(\xi_n/\|\xi_n\|_2)$ is the sequence in $l^2(\mathfrak{X})_{+}$ obtained from the sequence $\xi_n = \sum_{i \in \mathbb{N}} \alpha_i | \tilde{\xi}_n| \circ s_i^{-1} $. We denote by $\Q\tilde{\mathcal{P}}(\mathfrak{X})_\omega$ the set of internal probability measures obtained in this way and by $\Q\tilde{\mathcal{P}}^\infty(\mathfrak{X})_\omega$ the internal probability measures obtained in this way which annihilate $c_0 (\mathfrak{X})$.
\end{defn}

\begin{lem}
Every internal measure in $\Q \tilde{\mathcal{P}}(\mathfrak{X})_\omega$ is quasi-invariant.
\end{lem}
\proof The proof is analogous to the proof of Lemma \ref{lemza20} (cf. also the proof of \cite{BaRa} Theorem 2.7). $\Box$

Let $\mathcal{Q} (\mathfrak{X}) = \mathbb{B}(l^2 (\mathfrak{X})/\mathbb{K}(l^2(\mathfrak{X}))$ be the Calkin algebra of $\mathfrak{X}$.
\begin{defn}
Let $\phi \in \mathcal{S}(C^* (G))$. We say that $\phi$ is an \textit{$\mathfrak{X}$-Calkin state} if it factors through the canonical $*$-homomorphism $C^*(G) \rightarrow \mathcal{Q} (\mathfrak{X})$.
\end{defn}

\begin{defn}
Let $C^*_\tau (G)$ be a $C^*$-algebraic completion of the group algebra $\mathbb{C}[G]$ which is a quotient of the full group $C^*$-algebra $C^*(G)$ and let $\phi$ be a map from $C^*(G)$ to a set. We say that $\phi$ is {\it $\tau$-continuous} if it factors through the given surjection $C^*(G) \rightarrow C^*_\tau (G)$.
\end{defn}

We record the following result, which follows from the proof of \cite{BaRa} Theorem 2.7.

\begin{prop}
\label{prop2.14}
Let $\phi$ be an $\mathfrak{X}$-Calkin state on $G$. There are a free ultrafilter $\omega \in \partial_\beta \mathbb{N}$, a probability measure $\mu \in \Q \tilde{\mathcal{P}}^\infty(\mathfrak{X})_\omega$ and a vector $\eta \in L^2 (\mu)$ such that 
\begin{equation*}
\phi (s) = \langle u^\mu_s \eta, \eta \rangle_\mu \qquad \forall s \in G,
\end{equation*}
where $s \mapsto u^\mu_s$ is the Koopman representation associated to $\mu$.
\end{prop}
\proof Every state on $\mathcal{Q}(\mathfrak{X})$ is of the form $a \mapsto \lim_{n \rightarrow \omega} \langle x \xi_n, \xi_n \rangle$ for some free ultrafilter $\omega \in \partial_\beta \mathbb{N}$ and some sequence $ (\xi_n)$ in $l^2(\mathfrak{X})_{1}$ converging weakly to zero. Consider such a presentation for the state $\phi$. Let $(\alpha_i)$ be a sequence of positive real numbers such that $\alpha_i >0$ for every $i$, $\sum_{i,j} \alpha_i \alpha_j < \infty$. As in the previous section we consider the sequence $\tilde{\xi} = (\tilde{\xi}_n)$ given by
\begin{equation*}
\tilde{\xi}_n := \sum_{i \in \mathbb{N}}\alpha_i |\xi_n| \circ s_i^{-1}
\end{equation*}
This is a uniformly bounded sequence in $l^2 (\mathfrak{X})_+$. Consider the associated Radon measure $\tilde{\mu}=\mu_{\tilde{\xi}}$ on $\sigma(l^\infty (\mathfrak{X})_\omega)$ and let $\mu= \mu_{\tilde{\xi} /{\| \tilde{\mu}\|^{1/2}}}$. Since the original sequence $(\xi_n)$ converges to zero weakly, $\mu \in \Q \tilde{\mathcal{P}}^\infty (\mathfrak{X})_\omega$. Using Lemma 2.5 and Lemma 2.6 of \cite{BaRa}, arguing as in the proof of \cite{BaRa} Theorem 2.7 we see that the Koopman representation on $L^2(\mu)$ is given, on the dense linear suspace of $L^2(\mu)$ generated by elements in $l^\infty (\mathfrak{X})_\omega$ of the form $(\sum_{i=0}^N \alpha_i f_n (\xi_n \circ s_i^{-1}) / \tilde{\xi}_n)^\bullet$, by $[u^\mu_s (\sum_{i=0}^N \alpha_i f_n (\xi_n \circ s_i^{-1}) / \tilde{\xi}_n)]= [\sum_{i=0}^N \alpha_i (f_n \circ s^{-1}) (\xi_n \circ (s_i^{-1}s^{-1})) / \tilde{\xi}_n]$ (where by $[\cdot]$ we mean classes in $L^2(\mu)$). Hence the $\mathfrak{X}$-Calkin state $\phi$ is represented by $\phi (s) = \lim_{n\rightarrow \omega} \langle s \xi_n, \xi_n\rangle_{l^2(\mathfrak{X})} = \langle u^\mu_s \eta, \eta \rangle_\mu$ for $\eta = [ (\| \tilde{\mu}\|^{1/2} \xi_n / \tilde{\xi}_n)]$. $\Box$

\begin{prop}
\label{prop2.15}
The Calkin representation of $G$ on $\mathfrak{X}$ is $\tau$-continuous for some group $C^*$-algebra $C^*_\tau (G)$ if and only if for every $\omega \in \partial_\beta \mathbb{N}$ and every $\mu \in \Q \tilde{\mathcal{P}}^\infty (\mathfrak{X})_\omega$ the associated Koopman representation is $\tau$-continuous.
\end{prop}
\proof Suppose that the Calkin representation of $G$ on $\mathfrak{X}$ is $\tau$-continuous and let $\mu=\mu_\xi \in \Q \tilde{\mathcal{P}}^\infty (\mathfrak{X})_\omega$, where $\xi$ is associated to sequences $(\tilde{\xi}_n)$ and $(\alpha_i)$ as in Definition \ref{def2.2}. By the proof of Proposition \ref{prop2.14} the vectors of the form $[\sum_{i=0}^N \alpha_i f_n (\tilde{\xi}_n \circ s_i^{-1}) / \xi_n]$, with $(f_n)^\bullet \in l^\infty (\mathfrak{X})_\omega$, are dense in $L^2(\mu)$. Given such a vector, write it as $[(\eta_n / \tilde{\xi}_n)]$, where $\eta_n =  \sum_{i=0}^N \alpha_i f_n (\xi_n \circ s_i^{-1})$ for every $n$. Since $\eta_n$ converges to zero weakly with respect to $\omega$, the state $s \mapsto \lim_{n \rightarrow \omega} \langle s \eta_n/\|\eta_n\|, \eta_n /\|\eta_n\| \rangle_{l^2(\mathfrak{X})}$ is $\tau$-continuous. As we observed in the proof of Proposition \ref{prop2.14} this state is realized as $s \mapsto \langle u^\mu_s [\eta_n/(\tilde{\xi}_n \|\eta_n\|)], [\eta_n /(\tilde{\xi}_n \|\eta_n\|)]\rangle_\mu$ and so for every $a \in \mathbb{C}(G)$ we have $\| u^\mu (a)[\eta_n /\tilde{\xi}_n]\|^2_\mu = \langle u^\mu(a)u^\mu(a)^* [\eta_n /\tilde{\xi}_n], [\eta_n /\tilde{\xi}_n]\rangle_\mu= \lim_{n \rightarrow \omega}\langle  (aa^*)\eta_n/\|\eta_n\|, \eta_n /\|\eta_n\| \rangle_{l^2(\mathfrak{X})} \cdot \lim_{n \rightarrow \omega} \| \eta_n\|^2_{l^2(\mathfrak{X})} \leq \|aa^*\|_\tau \|[\eta_n /\tilde{\xi}_n]\|^2_\mu$. Hence the Koopman representation is $\tau$-continuous. On the other hand, if the Koopman representation is $\tau$-continuous for every $\mu \in \Q \tilde{\mathcal{P}}^\infty (\mathfrak{X})_\omega$, it follows that the Koopman representations constructed from $\mathfrak{X}$-Calkin states are $\tau$-continuous and so these states are $\tau$-continuous. $\Box$

\section{Statement of the result}

Before stating the main result, we note the following

\begin{lem}
\label{lemfinal}
Let $X$ be a compact Hausdorff space endowed with an action of a discrete countable group $G$. Let $\mu$ be a finite quasi-invariant Radon measure on $X$ which is the uniform limit of a net of finite uniformly bounded quasi-invariant Radon measures $\mu_\lambda$ which are Zimmer-amenable. Then $\mu$ is Zimmer-amenable.
\end{lem}
\proof We first observe the following general fact. Suppose that $\mu_\lambda$ is a net of uniformly bounded Radon measures converging uniformly (with respect to continuous functions) to a finite Radon measure $\mu$. We claim that then $\mu_\lambda$ converges to $\mu$ uniformly on uniformly bounded sets of Borel functions: $\sup_{h \; Borel, \|h\|_\infty \leq M} | \mu_\lambda (h) - \mu(h)| \rightarrow 0$ for every $M>0$. In fact, by regularity of the measures, the convergence is uniform on Borel sets. Hence it is uniform in the classes of simple functions $S_{R, N}:=\{ \sum_{i=1}^N \alpha_i \chi_{E_i}, \; E_i \mbox{ Borel set }, |\alpha_i|\leq R \; \forall i\}$. For any given $\epsilon >0$ and $M>0$ we can find a set $S_{R,N}$ such that for every $h$ Borel function with $\|h\|_\infty \leq R$ there is $s_h \in S_{R,N}$ such that $\| h-s\|_\infty < \epsilon$. Hence for every $\epsilon >0$ and $h$ Borel function with $\|h\|_\infty \leq M$ we have $|\int h d\mu_\lambda - \int h d\mu| \leq |\int h  d\mu_\lambda - \int s_h d\mu_\lambda| + |\int s_h d\mu_\lambda - \int s_h d\mu| + |\int s_h d\mu - h d\mu| \leq \sup_\lambda \mu_\lambda (X) \epsilon + \mu (X) \epsilon + |\int s_h d\mu_\lambda - \int s_h d\mu|$ and the result follows.


In virtue of \cite{BeCr} Theorem 3.6 it is sufficient to show that for every finite subset $A \subset L^\infty (X,\mu)_{*,+}$, every finite subset $F \subset G$ and every $\epsilon >0$ there is $g \in C_c(G, L^\infty (X,\mu)_+)$ such that
\begin{equation*}
\sum_{t \in G} g(t) =1, \qquad \omega (\sum_{t \in G} | sg (s^{-1} t) - g(t)|) < \epsilon \quad \forall \omega \in A, \; \forall s \in F.
\end{equation*}

Every $\omega \in L^\infty (X,\mu)_{*,+}$ can be represented by an element of $L^1(X,\mu)$, which we still denote by $\omega$. For every $\omega \in A$ choose a positive Borel representative $\tilde{\omega}$ and consider the Borel sets $E_{n,\omega} := \{ x \in X \; | \; \tilde{\omega} (x) \in [n,n+1)\}$. For every $\omega \in A$ there is $N_{\omega, \epsilon}$ such that $\int_{\bigsqcup_{n >N_{\omega , \epsilon}}E_{n, \omega}} \tilde{\omega} (x) d\mu < \epsilon/6$. Let $\tilde{\omega}_\epsilon:= \tilde{\omega} \chi_{\bigsqcup_{i=0}^{N-1} E_{i,\omega}}$. Then $\int | \tilde{\omega} - \tilde{\omega}_\epsilon| d\mu < \epsilon/6$ for every $\omega \in A$. Since $\mu_\lambda$ is a finite measure for every $\lambda$, the Borel functions $\tilde{\omega}_\epsilon$ are integrable for every $\lambda$. Hence we can find $\lambda$ such that there is $g_\lambda \in C_c (G, L^\infty (X,\mu_\lambda))_+$ satisfying
\begin{equation*}
\sum_{t \in G} g_\lambda(t) =1, \qquad \int_X \tilde{\omega}_\epsilon (\sum_{t \in G} | sg_\lambda (s^{-1} t) - g_\lambda(t)|) d\mu_\lambda< \frac{\epsilon}{3} \quad \forall \omega \in A, \; \forall s \in F.
\end{equation*}
For every $t$ in the support of $g_\lambda$ choose a positive Borel representative $\tilde{g}_\lambda (t)$ of $g_\lambda (t)$. We thus have $\sum_{t \in G} (\tilde{g}_\lambda (t) ) (x)=1$ for $\mu_\lambda$-almost every $x \in X$ and $ \int_X \tilde{\omega}_\epsilon (\sum_{t \in G} | s\tilde{g}_\lambda (s^{-1} t) - \tilde{g}_\lambda(t)|) d\mu_\lambda < \frac{\epsilon}{3} $ for every $\omega \in A$, $s \in F$. Let $E \subset X$ be the $\mu_\lambda$-null Borel set such that $\sum_{t \in G} (\tilde{g}_\lambda (t) ) (x) \neq 1$ and define $\hat{g}_\lambda(t) := \tilde{g}_\lambda (t)\chi_{X \backslash E} + |\supp (g_\lambda)|^{-1} \chi_E$. Then $\sum_{t \in G} g(t) =1$ for every $t$. Up to taking classes, we can consider $g$ as an element of $C_c (G, L^\infty (\mu))_+$ and $\tilde{\omega}_\epsilon$ as an element of both $L^1(X,\mu)$ and $L^1(X,\mu_\lambda)$. We have
\begin{equation*}
\begin{split}
\omega(\sum_{t \in G}& | s \hat{g}_\lambda(s^{-1} t) - \hat{g}_\lambda(t)|) \leq \|\omega - \tilde{\omega}_\epsilon\|_{L^1 (X,\mu)} (\| \sum_{t \in G} s \hat{g}_\lambda(s^{-1} t)\|_{L^\infty (X,\mu)} + \| \sum_{t \in G} \hat{g}_\lambda(t) \|_{L^\infty (X,\mu)})\\
& +\int \tilde{\omega}_\epsilon \sum_{t \in G} | s \hat{g}_\lambda(s^{-1} t) - \hat{g}_\lambda(t)| d\mu\\
&=  \|\omega - \tilde{\omega}_\epsilon\|_{L^1 (X,\mu)} (\| \sum_{t \in G} \hat{g}_\lambda( t)\|_{L^\infty (X,\mu)} + \| \sum_{t \in G} \hat{g}_\lambda(t) \|_{L^\infty (X,\mu)})
+\int \tilde{\omega}_\epsilon \sum_{t \in G} | s \hat{g}_\lambda(s^{-1} t) - \hat{g}_\lambda(t)| d\mu\\
& \leq 2 \|\omega - \tilde{\omega}_\epsilon\|_{L^1 (X,\mu)} + \int \tilde{\omega}_\epsilon \sum_{t \in G} | s \hat{g}_\lambda(s^{-1} t) - \hat{g}_\lambda(t)| d\mu\\
& \leq  2 \|\omega - \tilde{\omega}_\epsilon\|_{L^1 (X,\mu)} + |\int \tilde{\omega}_\epsilon \sum_{t \in G} | s \hat{g}_\lambda(s^{-1} t) - \hat{g}_\lambda(t)| d\mu - 
\int \tilde{\omega}_\epsilon \sum_{t \in G} | s \hat{g}_\lambda(s^{-1} t) - \hat{g}_\lambda(t)| d\mu_\lambda| \\ &+ \int \tilde{\omega}_\epsilon \sum_{t \in G} | s \hat{g}_\lambda(s^{-1} t) - \hat{g}_\lambda(t)| d\mu_\lambda\\ & < \frac{\epsilon}{3} +  |\int \tilde{\omega}_\epsilon \sum_{t \in G} | s \hat{g}_\lambda(s^{-1} t) - \hat{g}_\lambda(t)| d\mu - 
\int \tilde{\omega}_\epsilon \sum_{t \in G} | s \hat{g}_\lambda(s^{-1} t) - \hat{g}_\lambda(t)| d\mu_\lambda| +\frac{\epsilon}{3} .
\end{split}
\end{equation*}

Since the Borel functions $\tilde{\omega}_\epsilon \sum_{t \in G} | s \hat{g}_\lambda(s^{-1} t) - \hat{g}_\lambda(t)| $ are uniformly bounded, it follows from the observation at the start of the proof that for $\lambda$ large enough we have
\begin{equation*}
|\int \tilde{\omega}_\epsilon \sum_{t \in G} | s \hat{g}_\lambda(s^{-1} t) - \hat{g}_\lambda(t)| d\mu - 
\int \tilde{\omega}_\epsilon \sum_{t \in G} | s \hat{g}_\lambda(s^{-1} t) - \hat{g}_\lambda(t)| d\mu_\lambda| < \frac{\epsilon}{3} .
\end{equation*}
We conclude that if we pick $\lambda$ large enough we obtain
\begin{equation*}
\omega(\sum_{t \in G}s | s \hat{g}_\lambda(s^{-1} t) - \hat{g}_\lambda(t)|) < \epsilon \qquad \forall s \in F, \forall \omega \in A. \qquad  \Box
\end{equation*}

For the definitions of the sets of probability measures $\Q \tilde{\mathcal{P}}^\infty(\mathfrak{X})_\omega$ and $\Q \overline{F}^\infty_\omega (\mathfrak{X})$, see Definition \ref{def2.2} and Definition \ref{def02.2} respectively. $\Q \mathcal{P} (\partial_{\beta, \omega} \mathfrak{X})$ is the set of quasi-invariant probability measures on the non-standard extension of the Stone-\v{C}ech boundary $\partial_{\beta,\omega} \mathfrak{X}$ (see the discussion at the beginning of Subsection \ref{sec2}); equivalently, they are the quasi-invariant measures on the spectrum of $l^\infty (\mathfrak{X})_\omega$ which annihilate $c_0 (\mathfrak{X})$ under the canonical diagonal embedding (see Lemma \ref{lem2.6}). Recall that if $\mu$ is a quasi-invariant measure, we denote by $s \mapsto u^\mu_s$ the associated Koopman representation.

\begin{thm}
Let $\mathfrak{X}$ be a a countable set endowed with the action of a discrete countable group $G$ and let $\tau$ be a topology associated to a $C^*$-completion of the group algebra. Denote by $\pi_{\mathcal{Q}(\mathfrak{X})}$ the canonical representation of $G$ in $\mathcal{Q}(\mathfrak{X})$. Then
\begin{equation*}
\mbox{$\pi_{\mathcal{Q} (\mathfrak{X})}$ is $\tau$-continuous} \Leftrightarrow u^\mu \mbox{ is $\tau$-continuous } \forall \mu \in \Q \tilde{\mathcal{P}}^\infty(\mathfrak{X})_\omega, \; \forall \omega \in \partial_\beta \mathbb{N}
\end{equation*}
and
\begin{equation*}
G \curvearrowright \partial_\beta (\mathfrak{X}) \mbox{ is top. amenable} \Leftrightarrow \exists \omega \in \partial_\beta \mathbb{N} \; : \; \mbox{ $\mu$ is Zimmer-amenable } \forall \mu \in \Q \mathcal{P} (\partial_{\beta, \omega} (\mathfrak{X})).
\end{equation*}
Moreover, if $\pi_{\mathcal{Q}(\mathfrak{X})}$ is tempered, then every $\mu \in \Q \overline{F}^\infty_\omega (\mathfrak{X})$ is Zimmer-amenable.
\end{thm}
\proof The first claim is Proposition \ref{prop2.15}. The second claim follows from Corollary \ref{cor1.6} and Proposition \ref{prop3}. For the last claim, observe that every quasi-invariant Radon probability measure on $\sigma (l^\infty (\mathfrak{X})_\omega)$ which is a convex combination of uniform approximants of non-standard points associated to sequences $(\alpha^{(m)}_i)$ satisfying $\alpha^{(m)}_i>0$ for every $i$, $m$ is pair-amenable by Proposition \ref{propza5}; if it annihilates $c_0 (\mathfrak{X})$, it is also tempered by assumption, hence by Theorem \ref{thm3.5} it is Zimmer-amenable. Note that every measure in $\widetilde{\mathrm{NS}}^\infty_{\omega} (\mathfrak{X})$ is the norm limit of nets of such measures. Hence the same is true for every measure in $\Q \overline{F}^\infty_\omega (\mathfrak{X})$. The result follows since a quasi-invariant Radon probability measure which is in the norm closure of a set of Zimmer-amenable quasi-invariant Radon probability measures is Zimmer-amenable by Lemma \ref{lemfinal}. $\Box$

\section{Acknowledgments}
The author thanks Prof. R. Conti and Prof. F. R{\u a}dulescu for their useful comments on a previous version of this work. He is indebted to Prof. F. R{\u a}dulescu for many enlightening discussions concerning the topic of this research. He also thanks Prof. S. Echterhoff for an interesting discussion they had about the equivalence between measure-wise amenability and topological amenability during the conference "Group Actions: Dynamics, Measure, Topology", held at the University of M{\" u}nster, M{\" u}nster (DE) (28/11/2022 -- 02/12/2022). This research is supported by the grant "Algebre di Operatori e Teoria Quantistica dei Campi", CUP: E83C18000100006 and the University of Rome "Tor Vergata" funding OAQM, CUP: E83C22001800005. The author also acknowledges the support of INdAM-GNAMPA.

\bibliographystyle{mscplain}
 \bibliography{biblio}

\baselineskip0pt

\end{document}